\documentclass[twocolumn]{autart}  
\usepackage{graphicx}          
\usepackage{amsmath}
\usepackage{amssymb}
\usepackage{mathtools}
\usepackage{bm}
\usepackage{bbm}
\usepackage{mathrsfs}
\usepackage[sort&compress,square,numbers]{natbib}
\usepackage{xcolor}
\usepackage[pdftex, bookmarks=true, colorlinks=true, linkcolor=black, urlcolor=black, citecolor=black, linktoc=section]{hyperref}

\usepackage[normalem]{ulem}

\usepackage{enumerate}

\usepackage{subcaption}

\newcommand{\onetagright}{\tagsleft@false}
\newcommand\figref{Fig.~\ref}

\DeclareMathOperator{\cov}{cov}
\DeclareMathOperator{\vect}{vec}
\DeclareMathOperator{\kernl}{Ker}
\DeclareMathOperator{\img}{Im}

\let\oldforall\forall
\renewcommand{\forall}{\oldforall \,}

\newtheorem{definition}{Definition}[section]
\newtheorem{theorem}{Theorem}[section]
\newtheorem{lemma}{Lemma}[section]
\newtheorem{remark}{Remark}[section]
\newtheorem{proposition}{Proposition}[section]

\newtheorem{problem}{Problem}
\newtheorem{assumption}{Assumption}

\let\oldpf\pf
\let\oldendpf\endpf
\def\pf{\begingroup \oldpf}
\def\endpf{\hfill~\qed \oldendpf \endgroup}

\newcommand{\mR}{\mathbb{R}}

\newcommand{\mE}{\mathbb{E}}

\newcommand{\mS}{\mathbb{S}}

\DeclareMathOperator*{\argmin}{arg \, min}

\DeclareMathOperator*{\trace}{tr}

\newcommand{\ones}{\mathbf{1}}

\newcommand{\rvx}{\bm{x}}
\newcommand{\rvy}{\bm{y}}
\newcommand{\rvu}{\bm{u}}
\newcommand{\rvv}{\bm{v}}

\newcommand{\rvcx}{\bm{X}}
\newcommand{\rvcy}{\bm{Y}}
\newcommand{\rvcv}{\bm{V}}

\newcommand{\Qtrue}{\bar{Q}}
\newcommand{\Strue}{\bar{S}}
\newcommand{\Rtrue}{\bar{R}}
\newcommand{\Ptrue}{\bar{P}}
\newcommand{\Ktrue}{\bar{K}}

\graphicspath{{Figures/}}

\edef\endfrontmatter{%
  \unexpanded\expandafter{\endfrontmatter}
  \noexpand\endNoHyper 
}

\begin{document}

\begin{frontmatter}
  \title{Inverse linear-quadratic discrete-time finite-horizon optimal control for indistinguishable homogeneous agents: a convex optimization approach\thanksref{footnoteinfo}} 

  \thanks[footnoteinfo]{The work of Han Zhang was partially supported by the National Natural Science Foundation of China (NSFC), grant 62103276, and the work of Axel Ringh was supported by the Knut and Alice Wallenberg Foundation, grant KAW~2018.0349.}

  \author[SJTU,KEY-LAB,CENTRE]{Han Zhang}\ead{zhanghan\_tc@sjtu.edu.cn},
  \author[CHALMERS_AND_GU]{Axel Ringh\thanksref{HKUST}}\ead{axelri@chalmers.se}
  
  \thanks[HKUST]{The majority of the work was done when the author was with the Department of Electronic and Computer Engineering, The Hong Kong University of Science and Technology, Clear Water Bay, Kowloon, Hong Kong, China.}

  \address[SJTU]{Department of Automation, School of Electronic Information and Electrical Engineering, Shanghai Jiao Tong University, Shanghai, China}
  \address[KEY-LAB]{Key Laboratory of System Control and Information Processing, Ministry of Education of China, Shanghai 200240, China} 
  \address[CENTRE]{Shanghai Engineering Research Center of Intelligent Control and Management, Shanghai 200240, China}
  \address[CHALMERS_AND_GU]{Department of Mathematical Sciences, Chalmers University of Technology and University of Gothenburg, 41296 Gothenburg, Sweden}  


\begin{abstract}                          
The inverse linear-quadratic optimal control problem is a system identification problem whose aim is to recover the quadratic cost function and hence the closed-loop system matrices based on observations of optimal trajectories.
In this paper, the discrete-time, finite-horizon case is considered, where the agents are also assumed to be homogeneous and indistinguishable.
The latter means that the agents all have the same dynamics and objective functions and the observations are in terms of ``snap shots'' of all agents at different time instants, but what is not known is ``which agent moved where" for consecutive observations.
This absence of linked optimal trajectories makes the problem challenging.
We first show that this problem is globally identifiable. 
Then, for the case of noiseless observations, we show that the true cost matrix, and hence the closed-loop system matrices, can be recovered as the unique global optimal solution to a convex optimization problem.
Next, for the case of noisy observations, we formulate an estimator as the unique global optimal solution to a modified convex optimization problem.
Moreover, the statistical consistency of this estimator is shown.
Finally, the performance of the proposed method is demonstrated by a number of numerical examples. 
\end{abstract}

\end{frontmatter}

\section{Introduction}
First proposed in \citep{kalman1964linear}, Inverse Optimal Control (IOC) is an inverse problem whose aim is to reconstruct the cost function and predict the closed-loop system's behaviour, based on knowledge of the underlying system dynamics and observations of the system.
The problem can be categorized as a system identification problem where the system is assumed to be governed by an optimal control model with known structure, and hence it is a so-called gray-box system identification problem \cite[p.~13]{ljung1999system}.
IOC problems are of great interest, not least due to the fact that many processes in nature have been observed to be optimal with respect to some criteria \cite{alexander1996optima}. ``Reverse-engineering'' the agents' objective function not only enables us to understand how their decisions are made, but also provides a way to predict and imitate their behaviours.
For example, the IOC problem has found applications in human motion analysis \citep{berret2011evidence, westermann2020inverse}, robot locomotion \citep{mombaur2010human}, robot manipulation \citep{jin2019inverse, menner2019constrained}, systems biology \citep{tsiantis2018optimality} and machine learning \citep{finn2016guided, kopf2017inverse}.

Nevertheless, for some scenarios with multiple agents, such as human crowds, bacteria, bird flocks, and schools of fish, the observation data is naturally collected by observing all the agents simultaneously in discrete time, for instance via video cameras.
Due to the similarity between the agents and the fact that population is often huge, it is often too expensive or simply impossible to track each individual in the group.
Under such conditions, the agents can be assumed to be homogeneous and ``indistinguishable".
This ``indistinguishablility" means that we can not tell ``which agent goes where" for consecutive observations. ``Reverse-engineering" the objective function for such agents is challenging due to the fact that:
\begin{enumerate}
\item The control input measurement is usually not available.
Moreover, the ``indistinguishable" characteristics means that the agents' optimal trajectories are not available. Therefore, the existing IOC methods, such as \citep{zhang2019inverse,  yu2021system, molloy2018finite, jin2019inverse, zhang2019inverseCDC, aswani2018inverse}, that minimizes the difference between the observed and expected optimal trajectories, or minimizes the violation of optimality conditions, can not be applied.
\item The structural identifiability under such ``indistinguishable" observations is not guaranteed.
In particular, given observed data for optimal ``indistinguishable" agents, there is a priori no guarantee that there aren't different objective functions or closed-loop systems that corresponds to it.
\item The observations are usually contaminated by noise due to limitations of experiments and measuring methods. The identification methods should be statistically consistent and robust to observation noise in order to provide an accurate estimate of the objective function or the closed-loop system.
\end{enumerate}

On the other hand, the Linear-Quadratic (LQ) optimal control formulation is one of the most commonly used optimal control methods in practice. Due to its simplicity, it is often used to approximate complex control problems and optimal behaviours \citep{fridovich2020efficient, toumi2020tractable, ji2019feedback}.
It is therefore not surprising that the inverse LQ optimal control problem (with distinguishable observations) has been studied in different settings and in various degrees of generality; see, e.g., 
\citep[Sec.~5.6]{anderson2007optimal}, \citep[Sec.~10.5]{boyd1994linear} for the continuous-time infinite-horizon case,
\citep{li2018convex, li2020continuous} for the continuous-time finite-horizon case,
\citep{priess2014solutions} for the discrete-time infinite-horizon case,  
and \citep{keshavarz2011imputing, zhang2019inverse, zhang2019inverseCDC, yu2021system} for the discrete-time finite-horizon case, respectively. 
Moreover, although the inverse LQ optimal control problem can be seen as a special case of the IOC problem for general nonlinear systems \citep{molloy2018finite,  keshavarz2011imputing, molloy2020online, pauwels2016linear, hatz2012estimating, rouot2017inverse}, the special structure of the LQ formulation allows for precise answers to structural identifiability, and in the finite time horizon case for the identification of a linear time-varying system from a limited amount of data.
Further, to the best of our knowledge, the IOC problem for ``indistinguishable" homogeneous agents has not been considered in existing literatures.

In this work, we consider the linear-quadratic discrete-time finite-horizon IOC problem  with indistinguishable observations. In particular, we assume that each agent is governed by the same discrete-time finite-horizon Linear Quadratic Regulator (LQR). Moreover, in this initial work on ``indistinguishable'' homogeneous agents, we also assume that there is no interaction between the agents. This is a simplifying assumption which we intend to relax in the future. In the current work, we focus on studying how to tackle this ``indistinguishability'' issue in the observations.
In particular, this means that different trajectories of agents (which all follow the same dynamics) are observed, but such trajectories are ``indistinguishable" in the sense that the matching between each state observation and the corresponding agent is not available.
 The goal is to develop an inverse LQ optimal control algorithm for an external observer that can be used to identify the homogeneous agents' common objective function using ``indistinguishable'' observations.
More precisely, the main contributions are:
\begin{enumerate}[(i)]
\item We show that the corresponding model structure is strictly globally identifiable.
\item In the case of exact measurements, we formulate a convex optimization problem and prove that the unique globally optimal solution is the quadratic cost term of interest.
\item In the case of noisy observations, we formulate an estimator of the sought quadratic cost term as the unique globally optimal solution to a modified convex optimization problem. Moreover, we also show that this estimator is (asymptotically) statistically consistent. The latter means that the estimate converges to the true parameter value as the number of agents tends to infinity.
\end{enumerate}

The article is organized as follows: in Section~\ref{section:problem_formulation_and_identifiability}, we formulate the problem, specify the model structure, and show that the model structure is strictly globally identifiable. Section~\ref{sec:noiseless} investigates the inverse problem in the absence of measurement noise, and we formulate a convex optimization problem whose unique global optimal solution coincides with the parameters of interest. In Section~\ref{sec:noisy}, we adapt the results of the previous section to the case when the observations have measurement noise. In particular, it is shown that the unique global optimal solution of the modified stochastic optimization problem is a statistically consistent estimator. Next, discussions regarding some general numerical difficulties with the discrete-time finite-horizon LQ IOC problem are included in Section~\ref{sec:ill-cond}, and in Section~\ref{sec:numerics}, we present the numerical results. Finally, the article is concluded in Section~\ref{sec:conclusion}.

\textit{Notation:} The following notation is utilized throughout the article: $\mS^n$ denotes the set of $n\times n$ symmetric matrices, and $\mathbb{S}^n_+$ denotes the set of symmetric $n \times n$ matrices that are positive semidefinite. For symmetric matrices, $G_1\succeq G_2$ denotes the Loewner partial order of $G_1$ and $G_2$, namely, $G_1-G_2\in\mS^n_+$. Moreover, $G_1\succ G_2$ denotes that $G_1-G_2$ is strictly positive definite; see, e.g., \cite[Sec.~7.7]{horn2013matrix}.
$\|\cdot\|_F$ denotes the Frobenius norm, and $\ones$ denotes an all-one vector with appropriate length.
Furthermore, we use \textit{\textbf{italic bold font}} to denote stochastic elements.
Finally, we denote $\prod_{k=1}^m A_k = A_m A_{m-1}\ldots A_1$.

\section{Problem formulation, model structure, and identifiability}\label{section:problem_formulation_and_identifiability}
Suppose that each observed agent $i$ ($i=1,\ldots,M$) is governed by the following discrete-time finite-horizon LQR:
\begin{subequations}\label{eq:forward_problem}
  \begin{align}
  \min_{x_{1:N}^i, u_{1:N-1}^i}
  & \;\; x_N^{iT} \Strue x_N^{i} + \sum_{t = 1}^{N-1} \left( x_t^{iT} \Qtrue x_t^i + u_t^{iT} \Rtrue u_t^{i} \right) \label{eq:forward_problem_cost} \\
  \text{s.t.}
  & \;\; x_{t+1}^i = Ax_t^i + Bu_t^i, \; t = 1, \ldots, N-1, \label{eq:forward_problem_dynamics} \\
  & \;\; x_1^i = \bar{x}^i, \label{eq:forward_problem_init_cond}
  \end{align}
\end{subequations}
where $A, \Qtrue, \Strue \in \mR^{n \times n}$, $B \in \mR^{n \times m}$, and $\Rtrue \in \mR^{m \times m}$.
The agents are assumed to be homogeneous, namely, they have the same dynamics as well as the same objective functions. This means that the difference between the agents is only their initial values. For the sake of simplicity, as stated in the Introduction, in this work we assume there is no interactions between the agents. We intend to relax this assumption in future work. We will also assume that $\Strue = \Qtrue$ and $\Rtrue = I$ throughout this paper.  It is further assumed that $A$ is invertible, that $B$ has full column rank, and that $(A,B)$ is controllable. The rationale for making the assumption that $A$ is invertible lies in the fact that discrete-time systems are often obtained by sampling of continuous-time systems. More precisely, if $\hat{A}$ is the system matrix of the continuous-time system and $\Delta t > 0$ is the sample period, then the system matrix $A=\exp(\hat{A}\Delta t)$ of the sampled discrete-time systems is always invertible \citep{zhang2019inverse}.

The optimal control input sequence $\bar{u}_{1:N-1}^{i}$ of \eqref{eq:forward_problem} is given by $\bar{u}_t^{i} = \Ktrue_t x_t^{i}$, $t=1,\ldots,N-1$, where
\begin{align*}
\Ktrue_t = -(B^T\Ptrue_{t+1}B+I)^{-1}B^T\Ptrue_{t+1}A,
\end{align*}
and $\Ptrue_{2:N}$ is the solution to the following discrete-time Riccati Equation (DRE):
\begin{align}
\Ptrue_t = & \; A^T\Ptrue_{t+1}A+\Qtrue- \nonumber \\
&\; A^T\Ptrue_{t+1}B(B^T\Ptrue_{t+1}B+I)^{-1}B^T\Ptrue_{t+1}A, \nonumber \\
& \; t=2:N-1 \nonumber \\
\Ptrue_N = & \; \Qtrue. \label{eq:DRE}
\end{align}
For the problem of LQ IOC for indistinguishable agents, which is considered in this paper, the goal is to recover the matrix $\Qtrue$ in the objective function given the knowledge of (possibly noisy) state observations. However, note that since the agents are indistinguishable, this means that the observations at different time instants are \emph{unpaired}.
More precisely, 
let $X_t := [x_t^1, \ldots, x_t^M]$ be the state of $M$ agents following the dynamics in \eqref{eq:forward_problem}, then the observations are modelled as
\begin{align*}
\left\{Y_t^\phi \right\}_{t=1}^N=\left\{\left[x_t^{\phi_t(1)},\ldots, x_t^{\phi_t(M)}\right]\right\}_{t=1}^N,
\end{align*}
where $\phi_t:\{1,\ldots,M\} \mapsto \{1,\ldots,M\}$, $t = 1, \ldots, N$, are unknown permutations (the superscript $\phi$ in $Y_t^{\phi}$ is used to emphasizes the fact that the observations are shuffled).
Nevertheless, $\phi_1(\cdot)$ can, without loss of generality, be assumed to be the identity mapping since the ordering of the initial states does not matter. Therefore we will henceforth restrict our attention to the set $\{ \phi_t \}_{t = 2}^N$.

\begin{problem}{\textbf{(IOC for indistinguishable LQ homogeneous agents)}}\label{pro:problem_ioc}
Given the unpaired state observations $\{Y_t^\phi\}_{t=1}^N$ of $M$ homogeneous agents that are all governed by \eqref{eq:forward_problem}, recover the corresponding parameter $\Qtrue$ in the objective function.
\end{problem}

Before we continue, note that in the formulation in Problem~\ref{pro:problem_ioc} we assume that we have data from the entire planning horizon $N$ of the fowrad problem \eqref{eq:forward_problem}. This means that (implicitly) we assume that the planning horizon $N$ is known. Next, we would like to discuss the identifiability of such problem.  According to the definition of identifiability in \citep[Def.~4.5, 4.6]{ljung1999system},  \emph{the identifiability is a property of the model structure $\mathcal{M}$ itself and have nothing to do with any concrete data}, where the model structure $\mathcal{M}$ is a parameterized collection of models that describes the relations between the input and the output signal of the system \citep{ljung2013convexity}.  Hence, as a pre-step, we need to first define the underlying model structure $\mathcal{M}$ of Problem \ref{pro:problem_ioc}.

In Problem \ref{pro:problem_ioc}, we see the initial values of the trajectories as the ``input signal", and $\{Y_t^\phi\}_{t=2}^N$ as the ``output signal".
Note that any permutation $\phi_t(\cdot)$ can be uniquely represented by a permutation matrix $\mathscr{P}_t \in \mathcal{P}$, where
\[
\mathcal{P} := \{ \mathscr{P} \in \{0,1\}^{M \times M} \mid \mathscr{P}\ones = \ones, \;  \mathscr{P}^T\ones = \ones\},
\] 
i.e., the set of $0$-$1$-matrices with exactly one element that is one in each row and in each column. Post-multiplying any matrix with such a permutation matrix results in a matrix with permuted columns. Therefore, $Y^\phi_t$ can be represented as
\begin{align*}
Y^\phi_t = X_t \bar{\mathscr{P}}_t,\quad t=2,\cdots,N,
\end{align*}
where $\bar{\mathscr{P}}_t$ is the true (unknown) permutation, and hence for the system output it holds that
\begin{align}\label{eq:data_model_noise_free}
Y_t^\phi = \left(\prod_{k=1}^{t-1}A_{cl}(k;\Qtrue) X_1\right) \bar{\mathscr{P}}_t.
\end{align} 
where $A_{cl}(k;\Qtrue)$ is the closed-loop system matrix at time instant $k$ that is generated by LQR \eqref{eq:forward_problem} using the cost matrix $\Qtrue$.
By using the property of the vectorization operator \cite[Lem.~4.3.1]{horn1994topics}, vectorizing the above equation we have that
\begin{align*}
\vect(Y_t^\phi) &= \left[ \bar{\mathscr{P}}_t^T\otimes \prod_{k=1}^{t-1}A_{cl}(k; \Qtrue)\right]\vect(X_1).
\end{align*}
This is a valid input-output relation for any $Q \succeq 0$ and $\{\mathscr{P}_t\}_{t=2}^N \subset \mathcal{P}$.
Thus, we have the following definition of the model structure $\mathcal{M}$ for Problem \ref{pro:problem_ioc}:
\begin{align*}
\mathcal{M}(Q,\{\mathscr{P}_t\}_{t=2}^N)&:=\left\{\mathcal{M}_t(Q,\mathscr{P}_t)\right\}_{t=2}^N\\
&=\{\mathscr{P}_t^T\otimes \prod_{k=1}^{t-1}A_{cl}(k;Q)\}_{t=2}^N.
\end{align*}
Next, we adopt the definition of identifiability in \citep[Def. ~4.6, 4.7]{ljung1999system}. More precisely in our case:

\begin{definition}[Identifiability]
$\mathcal{M}$ is globally identifiable at $(Q^\prime,\{\mathscr{P}_t^\prime\}_{t=2}^N) \in \mathbb{S}^n_+ \times \mathscr{P}_t^{N-1}$ if it holds that 
\begin{equation*}
\begin{aligned}
&\mathcal{M}(Q,\{\mathscr{P}_t\}_{t=2}^N) = \mathcal{M}(Q^\prime,\{\mathscr{P}_t^\prime\}_{t=2}^N),\; Q\in\mathbb{S}^n_+, \mathscr{P}_t\in\mathcal{P}\\
&\implies Q=Q^\prime, \mathscr{P}_t = \mathscr{P}_t^\prime, \; \forall t=2,\ldots,N.
\end{aligned}
\end{equation*}
$\mathcal{M}$ is strictly globally identifiable if it is globally identifiable at all $Q^\prime\in\mathbb{S}^n_+,\mathscr{P}_t^\prime\in\mathcal{P}$.
\end{definition}

\begin{proposition}\label{prop:global_identifiability}
If $N \geq n +1$, then the model structure $\mathcal{M}(Q,\{\mathscr{P}_t\}_{t=2}^N)$ is strictly globally identifiable.
\end{proposition}
\begin{pf}
Let $(Q^\prime, \{\mathscr{P}_t^\prime\}_{t=2}^N) \in \mathbb{S}^n_+ \times \mathcal{P}^{N-1}$, and assume that $\mathcal{M}(Q,\{\mathscr{P}_t\}_{t=2}^N) = \mathcal{M}(Q^\prime,\{\mathscr{P}_t^\prime\}_{t=2}^N)$ for some $(Q, \{\mathscr{P}_t\}_{t=2}^N) \in \mathbb{S}^n_+ \times \mathcal{P}^{N-1}$.
Therefore, $\mathcal{M}_t(Q,\mathscr{P}_t)=\mathcal{M}_t(Q^\prime,\mathscr{P}_t^\prime)$ for $t=2,\ldots,N$. On the other hand, recall that any permutation matrix $\mathscr{P}_t$ is a 0-1 matrix with exactly one element that is one in each row and column. This means that for all $t$, $\mathcal{M}_t(Q,\mathscr{P}_t)$ is composed of all-zero and $\prod_{k=1}^{t-1}A_{cl}(k;Q)$ sub-matrix blocks. Together with the fact that $\mathcal{M}_t(Q,\mathscr{P}_t)=\mathcal{M}_t(Q^\prime,\mathscr{P}_t^\prime)$ for $t=2,\ldots,N$, it implies that $\mathscr{P}_t=\mathscr{P}_t^\prime$ and $\prod_{k=1}^{t-1} A_{cl}(k;Q)=\prod_{k=1}^{t-1} A_{cl}(k;Q^\prime)$ holds for $t=2,\ldots,N$. Since $A_{cl}(t;Q)$ is invertible for all $t=1,\ldots,N-1$ \citep{zhang2019inverse}, by induction $A_{cl}(t;Q)=A_{cl}(t;Q^\prime)$ holds for all $t=1,\ldots,N-1$. Now, if $N \geq n+1$, this implies that $Q=Q^\prime$ \citep[Theorem. ~2.1]{zhang2019inverse}. This gives the global identifiability of $\mathcal{M}(Q^\prime,\{\mathscr{P}_t^\prime\})$.
Finally, since $(Q^\prime, \{ \mathscr{P}_t^\prime\}_{t = 2}^N)$ was arbitrarily chosen from $\mathbb{S}^n_+ \times \mathcal{P}^{N-1}$,  $\mathcal{M}(Q^\prime,\{\mathscr{P}_t^\prime\}_{t=2}^N)$ is strictly globally identifiable.
\end{pf}

As a final note in this section, we observe that by identifying $\Qtrue$ we implicitly also identify $\{\bar{\mathscr{P}}_t\}_{t=2}^N$.
More precisely, if $\Qtrue$ is identified, then by solving the forward problem \eqref{eq:forward_problem}, the estimates of each agent's trajectory can be easily obtained; pairing these estimates with the data gives the permutations.
Moreover, as will be shown next, it turns out that $\Qtrue$ can be identified without explicitly identifying the permutations.

\section{IOC for indistinguishable agents in the noiseless case}
\label{sec:noiseless}

After justifying the identifiability of the problem, we now investigate the IOC problem for indistinguishable homogeneous agents in the noiseless case, namely, it is assumed that $y_t^i = x_t^{\phi_t(i)}$, $\forall t$ and $i$.
More precisely, we construct the IOC algorithm for indistinguishable observations as a convex optimization problem. However, before presenting the optimization problem, let us first sketch the main intuition behind it.

To this end, we first note that the optimal solution to the forward problem \eqref{eq:forward_problem} is characterized by the DRE \eqref{eq:DRE}. However, the latter is a nonlinear equation in $\Ptrue_{t}$, and in order to tackle to problem we first relax it to a linear matrix inequality (LMI).
To do so, note that since $\Qtrue \in \mathbb{S}_+^n$ we know that $\{ \Ptrue_t \}_{t=1}^N \subset \mathbb{S}^n_+$, and hence it follows that $B^T \Ptrue_t B+I\succ 0$, $\forall t=1:N$. Moreover, the expression $A^T \Ptrue_{t+1}A+ \Qtrue - A^T \Ptrue_{t+1}B(B^T \Ptrue_{t+1} B+I)^{-1}B^T \Ptrue_{t+1} A- \Ptrue_t$ is actually the Schur complement of 
\begin{equation*}
\bar{F}_t:=
\begin{bmatrix}
B^T \Ptrue_{t+1}B+I &B^T \Ptrue_{t+1}A\\
A^T \Ptrue_{t+1}B &A^T\Ptrue_{t+1}A+\Qtrue-\Ptrue_t
\end{bmatrix},
\end{equation*}
namely, $F_t\backslash(B^T \Ptrue_{t+1}B+I)$, which is well-defined since $B^T\Ptrue_{t+1}B+I$ has full rank.
By properties of the Schur complement (see, e.g., \cite[p.~495]{horn2013matrix}), if we let $F_t\succeq 0$ for $t=1,\ldots,N-1$, this is equivalent to relaxing the DRE into the matrix inequality
\begin{align}
&A^T \Ptrue_{t+1}A-\Ptrue_t+\Qtrue-A^T \Ptrue_{t+1}B(B^T\Ptrue_{t+1}B+I)^{-1}\nonumber\\
&\times B^T\Ptrue_{t+1}A\succeq 0.\label{eq:schur_complement_LMI_true}
\end{align}
On the other hand, introducing $\Ktrue_t:=-(B^T \Ptrue_{t+1}B+I)^{-1}B^T\Ptrue_{t+1}A$, the above relaxation can be further written as
\begin{align*}
(A+B \Ktrue_t)^T\Ptrue_{t+1}(A+B\Ktrue_t)-\bar{P}_t+\Qtrue+\Ktrue_t^T\Ktrue_t\succeq 0.
\end{align*}
Now, pre- and post-multiply on both sides of the above inequality with the state vector $x_t^i$, we get that
\begin{align*}
x_{t+1}^{iT}\Ptrue_{t+1}x_{t+1}^i-x_t^{iT}\Ptrue_t x_t^i+x_t^{iT}\Qtrue x_t^i+\|\bar{u}_t^i\|^2\ge 0,
\end{align*}
for $t=1,\ldots,N-1$, since $x_{t+1}^i = (A+B\Ktrue_t)x_t$ and $\bar{u}_t^i = \Ktrue_tx_t^i$.
Summing the above inequality from $t=1$ to $N-1$, and using that the trace operator is invariant under cyclic permutation of the elements, we get that
\begin{align*}
- \! \trace(\Ptrue_1x_1^ix_1^{iT}) \! + \! \trace(\Ptrue_Nx_N^ix_N^{iT}) \! + \! \! \sum_{t=1}^{N-1} \! \trace(\Qtrue x_t^ix_t^{iT}) \! + \! \|\bar{u}_t^i\|^2 \! \ge \! 0.
\end{align*}
Summing this inequality over $i$,  and using the property that
\begin{equation}\label{eq:YYT_is_XXT}
Y_{t}^{\phi} Y_{t}^{\phi T} = \sum_{i=1}^M x_{t}^{\phi_t(i)} x_{t}^{\phi_t(i) T} = \sum_{i=1}^M x_{t}^i x_{t}^{i T} = X_t^{} X_t^T,
\end{equation}
we get that
\begin{align*}
& -\trace(\Ptrue_1 Y_1^\phi Y_1^{\phi T})+\trace(\Ptrue_NY_N^\phi Y_N^{\phi T}) \\
& + \sum_{t=1}^{N-1} \Big( \trace(\Qtrue Y_t Y_t^{\phi T})+\sum_{i=1}^M\|\bar{u}_t^i\|^2  \Big) \ge 0.
\end{align*}
The above inequality is only derived for the ``true'' parameters $\Qtrue$ and $\{ \Ptrue_t \}_{t = 1}^N$ (in which case it is in fact an equality since \eqref{eq:schur_complement_LMI_true} is an equality in this case), but as we shall see shortly (Lemma~\ref{lem:H_lower_bound}) the inequality is in fact true for all parameters $Q \in \mathbb{S}_+^n$  and $\{ P_t \}_{t = 1}^N \subset \mathbb{S}_+^n$ such that
\begin{subequations}\label{eq:feasible_constraints}
  \begin{align}
  F_t :=&
    \begin{bmatrix}
    B^TP_{t+1}B+I &B^TP_{t+1}A\\
    A^TP_{t+1}B &A^TP_{t+1}A+Q-P_t
    \end{bmatrix}\succeq 0,\nonumber\\
    & t=1,\ldots,N-1,\label{eq:LMI}\\
    P_N = & Q.\label{eq:P_N_constr}
  \end{align}
\end{subequations}
Therefore, let us define the domain
\begin{equation}\label{eq:domain_D}
\mathscr{D} := \{Q \in \mS^n_+, \{P_t\in\mS^n_+\} : \eqref{eq:feasible_constraints} \text{ holds} \},
\end{equation}
and the function $H : \mathscr{D} \mapsto \mR$ as
\begin{align}\label{eq:def_H}
& H(Q,\{P_t\}) :=  -\trace(P_{1}Y_1^\phi Y_1^{\phi T}) \nonumber \\
& \: + \trace(P_NY_N^\phi Y_N^{\phi T}) +\sum_{t=1}^{N-1} \trace (QY_t^\phi Y_t^{\phi T}).
\end{align}
For this function, we have the following properties:

\begin{lemma}\label{lem:H_lower_bound}
For any $(Q,\{P_t\}_{t=1}^N) \in \mathscr{D}$, it holds that
\begin{align*}
H(Q,\{P_t\})+\sum_{i=1}^M\sum_{t=1}^{N-1}\|\bar{u}_t^{i}\|^2 \ge 0,
\end{align*}
where $\{\bar{u}_t^{i}\}_{t=1}^{N-1}$ is the optimal control input sequence in \eqref{eq:forward_problem} that corresponds to the ``true" matrix $\bar{Q}$ and the trajectory $\{x_t^{i}\}_{t=1}^N$.
Moreover, let  $\left\{\Ptrue_t\right\}_{t=1}^N$ be the solution to the DRE \eqref{eq:DRE} that corresponds to $\Qtrue$. Then $H(\Qtrue,\{\Ptrue_t\}) + \sum_{i=1}^M \sum_{t=1}^{N-1} \|\bar{u}_t^{i}\|^2 = 0$.
\end{lemma}

\begin{pf}
Since $\{P_t\}_{t=1}^N$, $Q$ is feasible, it is clear that $P_t\in\mathbb{S}^n_+$, $\forall t=1,\ldots,N$ and hence $B^TP_{t+1}B+I\succ 0,\forall t=1,\ldots,N-1$. Moreover, since \eqref{eq:LMI} holds, by positive semidefiniteness of the Schur complement \cite[p.~495]{horn2013matrix} it holds that
\begin{align}
&A^TP_{t+1}A-P_t+Q-A^TP_{t+1}B(B^TP_{t+1}B+I)^{-1}\nonumber\\
&\times B^TP_{t+1}A\succeq 0.\label{eq:schur_complement_LMI}
\end{align}
Introducing $K_t := -(B^TP_{t+1}B+I)^{-1}B^TP_{t+1}A$, we can re-write the above matrix inequalities as
\begin{align*}
(A+BK_t)^TP_{t+1}(A+BK_t)-P_t+Q+K_t^TK_t\succeq 0,\\
t=1,\ldots,N-1.
\end{align*}
Rearranging the terms,  we have
\begin{align}
&A^TP_{t+1}A-P_t+Q \succeq -K_t^TB^TP_{t+1}A \nonumber\\
&\quad -A^TP_{t+1}BK_t-K_t^T(B^TP_{t+1}B+I)K_t.\label{eq:expanded_LMI}
\end{align}
On the other hand, for the $i$'th trajectory, consider the following term
\begin{equation}
J_t^{i} = x_{t+1}^{i T}P_{t+1}x_{t+1}^{i} - x_t^{i T} P_t x_t^{i} + x_t^{i T} Q x_t^{i} + \|\bar{u}_t^{i}\|^2,\label{eq:J_def}
\end{equation}
where $\bar{u}_t^{i}$ is the control input that corresponds to $\bar{Q}$ and the trajectory $\{x_t^{i}\}_{t=1}^N$ in the forward problem \eqref{eq:forward_problem}. Using the fact that $x_{t+1}^{i}=(A+B\bar{K}_t)x_t^{i}$ and $\bar{u}_t^{i}=\bar{K}_t x_t^{i}$, it holds for \eqref{eq:J_def} that
\begin{align*}
&J_t^{i} =\nonumber\\
& x_t^{i T}\underbrace{\left[(A+B\bar{K}_t)^TP_{t+1}(A+B\bar{K}_t)-P_t+Q+\bar{K}_t^T\bar{K}_t\right]}_{=: \mathscr{H}_t}x_t^{i}.
\end{align*}
Expanding the terms of $\mathscr{H}_t$, we have
\begin{align*}
&\mathscr{H}_t = A^TP_{t+1}A-P_t+Q+\\
&\bar{K}_t^TB^TP_{t+1}A+A^TP_{t+1}B\bar{K}_t+\bar{K}_t^T(B^TP_{t+1}B+I)\bar{K}_t,
\end{align*}
and using the matrix inequality \eqref{eq:expanded_LMI} we get that 
\begin{align}
&\mathscr{H}_t \succeq \label{eq:LMI_H_t}\\
&-K_t^TB^TP_{t+1}A-A^TP_{t+1}BK_t-K_t^T(B^TP_{t+1}B+I)K_t \nonumber \\
&+\bar{K}_t^TB^TP_{t+1}A+A^TP_{t+1}B\bar{K}_t+\bar{K}^T(B^TP_{t+1}B+I)\bar{K}_t. \nonumber
\end{align}
Recalling that $K_t = -(B^TP_{t+1}B+I)^{-1}B^TP_{t+1}A$ and using the fact that
\begin{align*}
&\bar{K}_t^TB^TP_{t+1}A\\
&=\bar{K}_t^T(B^TP_{t+1}B+I)(B^TP_{t+1}B+I)^{-1}B^TP_{t+1}A\\
&=-\bar{K}_t^T(B^TP_{t+1}B+I)K_t,
\end{align*}
and that
\begin{align*}
&K_t^TB^TP_{t+1}A\\
&=K_t^T(B^TP_{t+1}B+I)(B^TP_{t+1}B+I)^{-1}B^TP_{t+1}A\\
&=-K_t^T(B^TP_{t+1}B+I)K_t,
\end{align*}
the matrix inequality \eqref{eq:LMI_H_t} can be further rewritten as
\begin{align*}
&\mathscr{H}_t\succeq 2K_t^T(B^TP_{t+1}B+I)K_t-K_t^T(B^TP_{t+1}B+I)K_t\\
&-\bar{K}_t^T(B^TP_{t+1}B+I)K_t-K_t^T(B^TP_{t+1}B+I)\bar{K}_t\\
&+\bar{K}_t^T(B^TP_{t+1}B+I)\bar{K}_t\\
&=(K_t-\bar{K}_t)^T(B^TP_{t+1}B+I)(K_t-\bar{K}_t).
\end{align*}
Hence it holds that
\begin{align}
&J_t^{i} = x_t^{i T}\mathscr{H}_t x_t^{i}\nonumber\\
&\ge x_t^{i T}(K_t-\bar{K}_t)^T(B^TP_{t+1}B+I)(K_t-\bar{K}_t)x_t^{i}.\label{eq:J_t_inequality}
\end{align}

On the other hand, recall that $Y_t^\phi=\left[x_t^{\phi_t(1)}, \cdots, x_t^{\phi_t(M)}\right]$. Therefore, it follows that
\begin{align}
&\trace(P_tY_t^\phi Y_t^{\phi T}) = \trace(P_t\left[x_t^{\phi_t(1)},\cdots, x_t^{\phi_t(M)}\right]\begin{bmatrix}
x_t^{\phi_t(1) T}\\ \vdots \\ x_t^{\phi_t(M) T}
\end{bmatrix}) \nonumber \\
&=\trace \left(P_t\left(\sum_{i=1}^M x_t^{\phi_t(i)} x_t^{\phi_t(i) T} \right)\right) = \sum_{i=1}^M x_t^{\phi_t(i) T}P_t x_t^{\phi_t(i)} \nonumber \\
&= \sum_{i=1}^M x_t^{i T}P_t x_t^{i},\quad t=1,\ldots,N. \label{eq:trace_P_Y_YT}
\end{align}
Similarly, it also holds that $\trace(Q Y_t^\phi Y_t^{\phi T})=\sum_{i=1}^M x_t^{i T}Q x_t^{i}$, $t=2,\ldots,N$. 
Together with \eqref{eq:J_def} and \eqref{eq:J_t_inequality}, we therefore have that
\begin{align}
& H(Q,\{P_t\}) + \sum_{i=1}^M \sum_{t=1}^{N-1} \|\bar{u}_t^{i}\|^2 = \sum_{i=1}^M \sum_{t=1}^{N-1} \|\bar{u}_t^{i}\|^2 + \nonumber \\
& \sum_{t=1}^{N-1} \left\{ \trace(P_{t+1}Y_{t+1}^\phi Y_{t+1}^{\phi T}) - \trace(P_tY_t^\phi Y_t^{\phi T}) + \trace(QY_t^\phi Y_t^{\phi T}) \right\} \nonumber \\
& = \! \sum_{t=1}^{N-1} \left\{ \sum_{i=1}^M \underbrace{x_{t+1}^{i T}P_{t+1}x_t^{i} - x_t^{i T} P_t x_t^{i} + x_t^{i T}Qx_t^{i} + \|\bar{u}_t^{i}\|^2}_{= J_t^{i}} \right\} \nonumber \\
& \ge \! \sum_{t=1}^{N-1} \! \left\{ \! \sum_{i=1}^M \left[x_t^{i T}(K_t-\bar{K}_t)^T(B^TP_{t+1}B+I)(K_t-\bar{K}_t)x_t^{i} \right] \! \right\} \nonumber \\
& \ge 0, \label{eq:H_plus_u_inequality}
\end{align}
which proves the first part of the statement.

Finally, by the arguments that lead up to \eqref{eq:def_H} we know that the lower bound is reached by $(\bar{Q},\left\{\bar{P}_t\right\}_{t=1}^N)$.
\end{pf}

In particular, note that $\{\bar{u}_t^{i}\}_{t=1}^{N-1}$ in Lemma~\ref{lem:H_lower_bound} has nothing to do with the variable $(Q,\{P_t\})$. Therefore, the lemma effectively shows that the function $H$ is bounded from below on the domain $\mathscr{D}$ and that the lower bound is attained by the ``true'' parameters $(\Qtrue, \left\{ \Ptrue_t \right\}_{t=1}^N)$. Therefore, we construct the following optimization problem to reconstruct $Q$:
\begin{subequations}\label{eq:opt_no_noise}
  \begin{align}
    & \underset{Q, \{ P_t \}_{t=1}^N }{\text{minimize}}
    & & H(Q,\{P_t\}) \label{eq:opt_no_noise_cost}\\
    & \text{subject to}
    & & (Q, \{ P_t \}_{t= 1}^N) \in \mathscr{D}.
  \end{align}
\end{subequations}
This is a semidefinite programming problem, and hence a convex optimization problem, that can be solved using standard numerical solvers, e.g., \citep{lofberg2004yalmip}.
We know that $\Qtrue$ and $\{ \Ptrue_t \}_{t = 1}^N$ is an optimal solution to the problem. Next, we show that it is unique. For this, we need the following result:

\begin{lemma}[Persistent excitation]\label{lem:persistent_excitation}
If $Y_1^\phi Y_1^{\phi T}\succ 0$, then it holds that $Y_t^\phi Y_t^{\phi T}\succ 0,\forall t=2,\ldots,N$.
\end{lemma} 
\begin{pf}
By using the property $x_{t+1}^{i} = A_{cl}(t; \Qtrue) x_t^{i}$, as well as \eqref{eq:YYT_is_XXT}, it follows that
\begin{align*}
Y_{t+1}^{\phi} Y_{t+1}^{\phi T} &= \sum_{i=1}^M x_{t+1}^{i} x_{t+1}^{i T}\\
&=A_{cl}(t; \Qtrue) \left(\sum_{i=1}^M x_t^{i} x_t^{i T}\right) A_{cl}^T(t; \Qtrue).
\end{align*}
Since for all $Q\in\mathbb{S}^n_+$, $A_{cl}(t;Q)$ is invertible for all $t=1,\ldots,N-1$ \citep{zhang2019inverse}, and since positive definiteness is invariant under congruence \cite[Obs.~7.1.8]{horn2013matrix}, 
by induction the statement follows.
\end{pf}

\begin{remark}
The ``persistent excitation" condition $Y_1^\phi Y_1^{\phi T}\succ 0$ is equivalent to that there are $n$ linearly independent initial values $x_1^i$ amongst all $M$ initial values.
\end{remark}

Now we are ready to present the main theorem for the noiseless case.

\begin{theorem}\label{thm:main_result_no_noise}
Assume that $N\ge n+1$, $Y_1^\phi Y_1^{\phi T}\succ 0$, and let $(Q^*,\{P_t^*\}_{t=1}^N)$ be an optimal solution to \eqref{eq:opt_no_noise}. Then $Q^*=\bar{Q}$, where $\bar{Q}$ is the ``true" parameter that is used to generate $\{Y_t^\phi\}_{t=2}^N$ with unknown permutations.
\end{theorem}
\begin{pf}
By Lemma \ref{lem:H_lower_bound}, we know that $\left(\bar{Q},\left\{\bar{P}_t\right\}_{t=1}^N\right)$ is an optimal solution to \eqref{eq:opt_no_noise}. What remains to show is the uniqueness.

Since $\left(Q^*,\{P_t^*\}_{t=1}^N\right)$ is an optimal solution, it must be feasible, thus \eqref{eq:H_plus_u_inequality} also holds for $\left(Q^*,\{P_t^*\}_{t=1}^N\right)$. Moreover, by the assumption that $Y_1^\phi Y_1^{\phi T}\succ 0$ and Lemma~\ref{lem:persistent_excitation}, it holds that $Y_t^\phi Y_t^{\phi T}\succ 0$, $\forall t=2,\ldots,N$, and hence there exists a strictly positive definite matrix $(Y_t^\phi Y_t^{\phi T})^{\frac{1}{2}}$ such that $Y_t^\phi Y_t^{\phi T}=(Y_t^\phi Y_t^{\phi T})^{\frac{1}{2}} (Y_t^\phi Y_t^{\phi T})^{\frac{1}{2}}$ \citep[p.~440]{horn2013matrix}.
Letting $K^*_t := -(B^TP_{t+1}^*B+I)^{-1}B^TP_{t+1}^*A$,
by further term manipulation in \eqref{eq:H_plus_u_inequality}, we have
\begin{align*}
& H(Q^*,\{P_t^*\})+\sum_{i=1}^M\sum_{t=1}^{N-1}\|\bar{u}_t^{i}\|^2\ge\\
& \sum_{t=1}^{N-1} \! \left\{ \! \sum_{i=1}^M \! \left[x_t^{i T}(K_t^*-\bar{K}_t)^T(B^TP_{t+1}^*B+I)(K_t^*-\bar{K}_t)x_t^{i}\right] \! \right\}\\
& = \! \sum_{t=1}^{N-1} \! \left\{ \! \trace \! \left[ (K_t^*-\bar{K}_t)^T(B^T P_{t+1}^*B+I)(K_t^*-\bar{K}_t)Y_t^\phi Y_t^{\phi T} \right] \! \right\} \\
& = \! \sum_{t=1}^{N-1} \! \left\{ \left\|(B^T P_{t+1}^*B+I)^{\frac{1}{2}}(K_t^*-\bar{K}_t)(Y_t^\phi Y_t^{\phi T})^{\frac{1}{2}}\right\|_F^2 \! \right\} \\
& \ge 0.
\end{align*}
As stated in Lemma \ref{lem:H_lower_bound},  the lower bound zero in the above inequality is reached by $\left(\bar{Q},\left\{\bar{P}_t\right\}_{t=1}^N\right)$, and it also holds that $H(Q^*,\{P_t^*\}_{t=1}^N)+\sum_{i=1}^M\sum_{t=1}^{N-1}\|\bar{u}_t^{i}\|^2=0$ since $\left(Q^*,\{P_t^*\}_{t=1}^N\right)$ is also an optimal solution. Therefore, since all terms in the above sum are nonnegative, it must hold that 
\begin{align*}
&\left\|(B^T P_{t+1}^*B+I)^{\frac{1}{2}}(K_t^*-\bar{K}_t)(Y_t^\phi Y_t^{\phi T})^{\frac{1}{2}}\right\|_F=0,\\
&\Leftrightarrow (B^T P_{t+1}^*B+I)^{\frac{1}{2}}(K_t^*-\bar{K}_t)(Y_t^\phi Y_t^{\phi T})^{\frac{1}{2}}=0,\\
&\Leftrightarrow K_t^*-\bar{K}_t=0,\quad t=1,\ldots,N-1
\end{align*}
since $(B^T P_{t+1}^*B+I)^{\frac{1}{2}}$ and $(Y_t^\phi Y_t^{\phi T})^{\frac{1}{2}}$ are invertible.

The above argument shows that $K_t^*=\bar{K}_t$ for all $t=1,\ldots,N-1$. Nevertheless, note that $\{K_t^*\}_{t=1}^{N-1}$ is not the control gain that is generated by $Q^*$ using DRE \eqref{eq:DRE}. Instead, it is defined by $K_t^* = -(B^TP_{t+1}^*B+I)^{-1}B^TP_{t+1}^*A$, for $t=1,\ldots,N-1$, where $(Q^*,\{P_t^*\}_{t=1}^{N-1})$ is an optimizer of \eqref{eq:opt_no_noise}. Hence the result from \citep[Theorem~ 2.1]{zhang2019inverse} can not be directly applied to show that $Q^*=\bar{Q}$. 
Nevertheless, to show the latter, first note that we can always write $P_t^* = \bar{P}_t+\Delta P_t,t=1,\ldots,N$, and $Q^* = \bar{Q}+\Delta Q$, for some symmetric matrices $\Delta Q$ and $\Delta P_t$, for $t = 1, \ldots, N$. Since both $(Q^*,\{P_t^*\}_{t=1}^N)$ and $(\bar{Q},\{\bar{P}_t\}_{t=1}^N)$ are optimal solutions to \eqref{eq:opt_no_noise}, they must have the same optimal value, and hence
\begin{align}
&H(Q^*,\{P_t^*\}_{t=1}^N) - H(\bar{Q},\{\bar{P}_t\}_{t=1}^N)=0\nonumber\\
\implies& -\trace(\Delta P_1 Y_1^\phi Y_t^{\phi T})+\trace(\Delta P_N Y_N^\phi Y_N^{\phi T}) \nonumber\\
&+ \sum_{t=1}^{N-1}\trace(\Delta Q Y_t^\phi Y_t^{\phi T})=0\label{eq:Delta_P_1_Delta_P_N_Delta_Q},
\end{align}
where we have used the linearity of the trace operator.
On the other hand, it holds that
\begin{align*}
& \quad K_t^*=-(B^TP_{t+1}^*B+I)^{-1}B^TP_{t+1}^*A\\
\Longleftrightarrow & \quad (B^TP_{t+1}^*B+I)K_t^* = -B^TP_{t+1}^*A\\
\Longleftrightarrow & \quad B^TP_{t+1}^*(A+BK_t^*) = -K_t^*.
\end{align*}
Since $K_t^* = \bar{K}_t,\forall t=1,\ldots,N-1$, $A+BK_t^*=A+B\bar{K}_t = A_{cl}(t; \bar{Q})$, and by the fact that $A_{cl}(t;\bar{Q})$ is invertible for all $t=1,\ldots, N-1$ \citep{zhang2019inverse}, it follows that
\begin{align}
&B^T P_{t+1}^* = -K_t^*A_{cl}(t;\bar{Q})^{-1} = -\bar{K}_tA_{cl}(t;\bar{Q})^{-1}=B^T\bar{P}_{t+1}\nonumber\\
&\implies B^T\Delta P_{t+1} = 0,\forall t=1,\ldots,N-1.\label{eq:B^T_Delta_P}
\end{align}

Next, since $(Q^*,\{P_t^*\}_{t=1}^N)$ is feasible, \eqref{eq:schur_complement_LMI} also holds for $(Q^*,\{P_t^*\}_{t=1}^N)$. This means that
\begin{align}
&A^T(\bar{P}_{t+1}+\Delta P_{t+1})A-(\bar{P}_t+\Delta P_t) + (\bar{Q}+\Delta Q)\nonumber\\
&-A^T(\bar{P}_{t+1}+\Delta P_{t+1})B[B^T(\bar{P}_{t+1}+\Delta P_{t+1})B+I]^{-1}\nonumber\\
&B^T(\bar{P}_{t+1}+\Delta P_{t+1})A\succeq 0,\:t=1,\ldots,N-1\nonumber\\
&\bar{P}_N+\Delta P_N = \bar{Q}+\Delta Q.\nonumber
\end{align}
On the other hand, since $(\bar{Q},\{\bar{P}_t\}_{t=1}^N)$ satisfies \eqref{eq:DRE}, by also using \eqref{eq:B^T_Delta_P} we get that
\begin{align}
&A^T\Delta P_{t+1}A-\Delta P_t + \Delta Q \succeq 0,\;\; t=1,\ldots,N-1,\label{eq:Delta_P_Delta_Q_inequality}\\
&\Delta P_N=\Delta Q.\label{eq:Delta_P_Delta_Q}
\end{align}
Using \eqref{eq:B^T_Delta_P}, we can further manipulate the terms of \eqref{eq:Delta_P_Delta_Q_inequality} and get
\begin{align}
&\underbrace{(A+B\bar{K}_t)^T}_{A_{cl}(t;\bar{Q})}\Delta P_{t+1} \underbrace{(A+B\bar{K}_t)}_{A_{cl}(t;\bar{Q})} - \Delta P_t + \Delta Q\succeq 0,\nonumber\\
&\qquad\forall t=1,\ldots,N-1\nonumber\\
\implies &X_t^TA_{cl}(t;\bar{Q})^T\Delta P_{t+1}A_{cl}(t;\bar{Q})X_t - X_t^T\Delta P_t X_t \nonumber\\
&+ X_t^T\Delta Q X_t\succeq 0,\nonumber\\
\implies &X_{t+1}^T\Delta P_{t+1}X_{t+1}-X_t^T\Delta P_t X_t + X_t^T\Delta Q X_t\succeq 0,\nonumber\\
\implies &\trace(\Delta P_{t+1}X_{t+1}X_{t+1}^T)-\trace(\Delta P_t X_tX_t^T)\nonumber\\
&+\trace(\Delta QX_tX_t^T)\ge 0,t=1,\ldots,N-1, \nonumber\\
\implies &\trace(\Delta P_{t+1}Y_{t+1}^\phi Y_{t+1}^{\phi T})-\trace(\Delta P_t Y_t^\phi Y_t^{\phi T})\nonumber\\
&+\trace(\Delta QY_t^\phi Y_t^{\phi T})\ge 0,t=1,\ldots,N-1, \label{eq:trace_delta_ineq}
\end{align}
where $X_t=[x_t^1,\cdots,x_t^M]$, and where we also used that $\trace(\Delta P_tY_t^\phi Y_t^{\phi T}) = \trace(\Delta P_tX_tX_t^T)$ and that $\trace(\Delta P_{t+1}Y_{t+1}^\phi Y_{t+1}^{\phi T}) = \trace(\Delta P_{t+1}X_{t+1}X_{t+1}^T)$.

Summing \eqref{eq:trace_delta_ineq} from $t=1$ to $N-2$,  we have
\begin{align}
&\trace(\Delta P_{N-1}Y_{N-1}^\phi Y_{N-1}^{\phi T})-\trace(\Delta P_1 Y_1^\phi Y_1^{\phi T})\nonumber\\
&\qquad+\sum_{t=1}^{N-2}\trace(\Delta QY_t^\phi Y_t^{\phi T})\ge 0. \label{eq:trace_delta_ineq_sum}
\end{align}
In view of \eqref{eq:Delta_P_1_Delta_P_N_Delta_Q} and \eqref{eq:trace_delta_ineq_sum}, we have that
\begin{align*}
0=&-\trace(\Delta P_1 Y_1^\phi Y_1^{\phi T})+\trace(\Delta P_N Y_N^\phi Y_N^{\phi T}) \\
&+ \sum_{t=1}^{N-1}\trace(\Delta Q Y_t^\phi Y_t^{\phi T})\\
=&-\trace(\Delta P_1 Y_1^\phi Y_1^{\phi T})+\trace(\Delta P_{N-1} Y_{N-1}^\phi Y_{N-1}^{\phi T})\\
&+ \sum_{t=1}^{N-2}\trace(\Delta Q Y_t^\phi Y_t^{\phi T})-\trace(\Delta P_{N-1}Y_{N-1}^\phi Y_{N-1}^{\phi T})\\
&+\trace(\Delta P_N Y_N^\phi Y_N^{\phi T})+\trace(\Delta QY_{N-1}^\phi Y_{N-1}^{\phi T})\\
\ge&-\trace(\Delta P_{N-1}Y_{N-1}^\phi Y_{N-1}^{\phi T})+\trace(\Delta P_N Y_N^\phi Y_N^{\phi T})\\
&+\trace(\Delta QY_{N-1}^\phi Y_{N-1}^{\phi T})\\
=&-\trace(\Delta P_{N-1}X_{N-1}X_{N-1}^T)+\trace(\Delta P_N X_N X_N^T)\\
&+\trace(\Delta QX_{N-1}X_{N-1}^T)
\end{align*}
Since $X_N=A_{cl}(N-1;\bar{Q})X_{N-1}$ and by also using \eqref{eq:B^T_Delta_P}, from the equation above it follows that
\begin{align*}
0\ge& \trace\Big[ \Big(A_{cl}(N-1;\bar{Q})^T\Delta P_N A_{cl}(N-1;\bar{Q})-\Delta P_{N-1}\\
&+\Delta Q \Big)X_{N-1}X_{N-1}^T\Big]\\
=&\trace\Big[(A^T\Delta P_N A-\Delta P_{N-1}+\Delta Q)X_{N-1}X_{N-1}^T\Big].
\end{align*}
Since $X_{N-1}X_{N-1}^T$ is positive definite (by Lemma~\ref{lem:persistent_excitation}), it follows that $A^T\Delta P_N A-\Delta P_{N-1}+\Delta Q\preceq 0$. Together with \eqref{eq:Delta_P_Delta_Q_inequality}, we can therefore conclude that $A^T\Delta P_N A-\Delta P_{N-1}+\Delta Q= 0$. 
Using \eqref{eq:B^T_Delta_P}, we therefore have that $A_{cl}(N-1; \bar{Q})^T\Delta P_N A_{cl}(N-1; \bar{Q})-\Delta P_{N-1}+\Delta Q= 0$, which implies that $-\trace(\Delta P_{N-1}X_{N-1}X_{N-1}^T)+\trace(\Delta P_N X_N X_N^T)+\trace(\Delta QX_{N-1}X_{N-1}^T)=0$. Therefore, summing \eqref{eq:trace_delta_ineq} from $t=1$ to $N-3$ and reiterating the above analysis, we can conclude that $A^T\Delta P_{N-1} A-\Delta P_{N-2}+\Delta Q= 0$. Doing this recursively, we have that
\begin{align}
A^T\Delta P_{t+1} A-\Delta P_t+\Delta Q= 0,\; t=1,\ldots, N-1.\label{eq:A^T_Delta_PA_Delta_P_Delta_Q}
\end{align}
Equipped with \eqref{eq:A^T_Delta_PA_Delta_P_Delta_Q}, we can use the same argument as in the proof of \citep[Theorem 2.1]{zhang2019inverse} to conclude that $\Delta Q=0$. Thus $Q^*=\bar{Q}$,  i.e., the optimal solution of \eqref{eq:opt_no_noise} is unique and hence the theorem follows.
\end{pf}

Before proceeding, we make a few remarks about the formulation in \eqref{eq:opt_no_noise}. As will be seen later, these remarks naturally hold also in the case of noisy observations.

\begin{remark}\label{rem:distinguishable_and_indistinguishable}
The optimization problem \eqref{eq:opt_no_noise} can also be applied to the IOC problem for ``distinguishable" agents, i.e., when the trajectories of the agents are available. Nevertheless, as one can see from \eqref{eq:opt_no_noise}, in this formulation of the problem there is no fundamental difference between IOC problems for ``distinguishable" and ``indistinguishable" agents.
\end{remark}

\begin{remark}\label{rem:scaling_num_agents}
The size of the variables and the number of variables and constraints in the optimization problem \eqref{eq:opt_no_noise} only depends on the size of the state space, $n$, and the length of the time horizon, $N$, respectively.
In particular, the overall size of the problem \eqref{eq:opt_no_noise} is independent of the number of agents, $M$. 
The only quantity that scales with $M$ is the computations of the outer products $Y_t^\phi Y_t^{\phi T}$, $t = 1,\ldots, N$. The latter scales linearly in $M$ and can be done before solving the optimization problem.
Therefore, the problem can be efficiently solved for scenarios with a large number of agents.
\end{remark}

\section{IOC for indistinguishable agents in the noisy case}\label{sec:noisy}
Next, we extend the results to the case of noisy observations of the unpaired states. In particular, we show that a modified version of \eqref{eq:opt_no_noise} gives a statistically consistent estimate.
To this end, let $(\Omega, \mathcal{F},\mathbb{P})$ be a probability space which carries random vectors $\rvx_1^i\in\mR^n$, $\{\rvv_t^i \in \mR^n\}_{t = 1}^N$, for  $i = 1,\ldots, M$. Based on that, the following assumptions are made for the rest of the paper.

\begin{assumption}[I.I.D.~random variables]\label{ass:noise_iid_and_zero_mean}
The random vectors $\rvx_1^i\in\mR^n$, $\{\rvv_t^i \in \mR^n\}_{t = 1}^N$, for  $i = 1,\ldots, M$, are all independent. Moreover,  random vectors  $\rvx_1^i\in\mR^n$, for  $i = 1,\ldots, M$, are identically distributed, and the random vectors $\rvv_t^i \in \mR^n$, for $t = 1, \ldots, N$ and $i = 1,\ldots, M$, are  identically distributed. Finally, $\mE[\rvv_t^{i}]=0$, and $\cov(\rvv_t^{i}, \rvv_t^{i}):=\Sigma$ is a priori known, where $\|\Sigma\|_F<\infty$.
\end{assumption}

\begin{assumption}[Stochastic persistent excitation]\label{ass:persistent_excitation_noisy}
It holds that $\cov(\rvx_1^{i},\rvx_1^{i})\succ 0$, where $\mE[\|\rvx_1^{i}\|^2]<\infty$.
\end{assumption}

\begin{assumption}[Bounded parameter]\label{ass:Q_set_compact}
The ``true" $\bar{Q}$ that governs the agents lives in a compact set $\bar{\mathbb{S}}^n_+(\varphi) := \{Q\in\mathbb{S}^n_+ : \|Q\|_F^2\le \varphi\}$, for some $0 < \varphi < \infty$.
\end{assumption}

Our goal is to seek the ``true" $\bar{Q}$ in $\bar{\mathbb{S}}^n_+(\varphi)$. To this end, we define a domain $\mathscr{D}(\varphi)$ similar to $\mathscr{D}$ in \eqref{eq:domain_D}, namely,
\begin{equation}\label{eq:domain_D_varphi}
\mathscr{D}(\varphi) := \{Q\in\bar{\mS}^n_+(\varphi),\{P_t\in\mS^n_+\}:\eqref{eq:LMI},\eqref{eq:P_N_constr}\}.
\end{equation}
Note that in practice, we can set $\varphi$ arbitrarily large if there is no prior knowledge on the norm bound of possible $Q$.

Equipped with the stochastic problem set-up, let the initial value $x_1^{i}$ be a realization of $\rvx_1^{i}$. Then the optimal state trajectory and optimal control input of the ``forward" optimal control \eqref{eq:forward_problem} can be seen as mappings from $\Omega$ to $\mR^n$ and $\mR^m$, respectively.
This means that the states and control signals are in fact stochastic variables, which are parameterized by $Q$ and implicitly determined by
\begin{align}
&\{\bar{\rvx}_t^{i}(\omega)\}_{t=2}^N,\{ \bar{\rvu}_t^{i}(\omega)\}_{t=1}^{N-1} = \label{eq:stochastic_forward_problem} \\
& \quad \left\{
\begin{array}{cl}
\argmin & \; J\left(\{\rvx_t^{i}(\omega)\}_{t=2}^N,\{\rvu_t^{i}(\omega)\}_{t=1}^{N-1}; \Qtrue, \rvx_1^{i}(\omega) \right) \nonumber \\
\text{subject to} & \; \eqref{eq:forward_problem_dynamics}, \eqref{eq:forward_problem_init_cond}, \nonumber
\end{array}
\right.
\end{align}
for all $\omega\in \Omega$, and where $J$ is the cost function in \eqref{eq:forward_problem_cost}.
From now on, we omit the ``bar'' and simply write $\{\rvx_t^{i}\}$ to denote the corresponding optimal state for the sake of brevity.
Next, we assume that the noisy observations of a state is given by  $\tilde{\rvx}_t^i = \rvx_t^i + \rvv_t^i$, and that the observed data thus take the form 
$\rvcy_t^{\phi} = [\rvy_t^{1}, \ldots, \rvy_t^{M}] = [\tilde{\rvx}_t^{\phi_t(1)}, \ldots, \tilde{\rvx}_t^{\phi_t(M)}]$.
Written in a form similar to \eqref{eq:data_model_noise_free}, the measured output can be expressed as
\begin{align*}
\rvcy_t^\phi = \left( \rvcx_t + \rvcv_t \right) \bar{\mathscr{P}}_t = \! \left( \! \left( \prod_{k=1}^{t-1}A_{cl}(k; \Qtrue) \rvcx_1 \! \right) \! + \! \rvcv_t \! \right) \bar{\mathscr{P}}_t,
\end{align*}
where $\rvcx_t = [\rvx_t^1, \ldots, \rvx_t^M]$ and $\rvcv_t = [\rvv_t^1, \ldots, \rvv_t^M]$.
This together with Assumption~\ref{ass:noise_iid_and_zero_mean} implies that the columns of $\rvcy_t^\phi$ are I.I.D.
A calculation similar to \eqref{eq:trace_P_Y_YT} then shows that
\begin{equation*}
\mE\left[ \trace \left( Q \rvcy_t^\phi \rvcy_t^{\phi T} \right)  \right] = \mE\left[ \trace \left( Q \rvcx_t \rvcx_t^{T} \right) \right] + M\trace \left( Q \Sigma \right)
\end{equation*}
where $\Sigma$ is the covariance of the noise.
In particular, the last equality holds since $\rvx_t ^{i}= \prod_{k = 1}^{t-1}A_{cl}(k; \Qtrue) \rvx_1^{i}$, and hence $\rvx_t^{i}$ is a stochastic variable that is independent of $\rvv_t^{i}$, for $t = 1, \ldots, N$ and $i = 1,\ldots M$.
Similar expressions hold for the other terms in \eqref{eq:def_H}, and based on this we construct the problem
\begin{subequations}\label{eq:opt_noise}
  \begin{align}
    & \underset{Q,  \{P_t \}_{t=1}^N}{\text{minimize}}
    & & \mE\left[ H_S^{\rvcy}(Q,\{P_t\}) \right] \label{eq_opt_noise_cost} \\
    & \text{subject to}
    & & (Q,  \{P_t \}_{t=1}^N) \in \mathscr{D}(\varphi),
  \end{align}
where
\begin{align}\label{eq:def_H_stochastic}
& H_S^{\rvcy}(Q,\{P_t\}) := \frac{1}{M}\Bigg[ -\trace(P_{1}\rvcy_1^{\phi } \rvcy_1^{\phi T}) \nonumber \\
& + \trace(P_N \rvcy_N^{\phi } \rvcy_N^{\phi T}) +\sum_{t=1}^{N-1} \trace (Q \rvcy_t^{\phi } \rvcy_t^{\phi T}) \nonumber \\
& + M\big( \trace \left( P_1\Sigma \right) - \trace(P_N \Sigma) - (N-1)\trace (Q\Sigma) \big) \Bigg].
\end{align}
\end{subequations}
In particular, note that by  a direct calculation it follows that the cost function \eqref{eq_opt_noise_cost} can be written as
\begin{align}
& \mE\left[ H_S^{\rvcy}(Q,\{P_t\}) \right] = \frac{1}{M}\mE\Big[-\trace(P_{1}\rvcx_1 \rvcx_1^{T}) \nonumber \\
& + \trace(P_N \rvcx_N \rvcx_N^{T}) + \sum_{t=1}^{N-1} \trace (Q \rvcx_t \rvcx_t^{T}) \Big].  \label{eq:H_S^y_equivalence}  
\end{align}
With this, we can now prove the following result akin Theorem~\ref{thm:main_result_no_noise}.

\begin{proposition}\label{prop:stoc_prob_unique_sol}
Let $\bar{Q}$ be the true parameter in \eqref{eq:stochastic_forward_problem}, and let $\left\{\bar{P}_t\right\}_{t=1}^N$ be the corresponding solution to the DRE \eqref{eq:DRE}. Under Assumptions~\ref{ass:noise_iid_and_zero_mean}, \ref{ass:persistent_excitation_noisy}, and \ref{ass:Q_set_compact}, if $N\ge n+1$, then $\left(\bar{Q},\left\{\bar{P}_t\right\}_{t=1}^N\right)$ is the unique optimal solution to \eqref{eq:opt_noise}.
\end{proposition}

\begin{pf}
Since \eqref{eq:H_S^y_equivalence} holds, by adapting the arguments in the proof of Lemma~\ref{lem:H_lower_bound} it follows that
\begin{align}
&\mE[H_S^{\rvcy}(Q,\{P_t\})] + \frac{1}{M}\sum_{i = 1}^M \sum_{t=1}^{N-1}\mE\left[\| \bar{\rvu}_t^{i}\|^2\right]\nonumber\\
&\ge \frac{1}{M}\sum_{i = 1}^M \sum_{t=1}^{N-1} \mE \left[\rvx_t^{iT}(K_t-\bar{K}_t)^T(B^TP_{t+1}B+I)(K_t-\bar{K}_t)\rvx_t^{i}\right]\nonumber\\
&=\frac{1}{M}\sum_{t=1}^{N-1} \trace\left[ (K_t-\bar{K}_t)^T(B^TP_{t+1}B+I)(K_t-\bar{K}_t)\mE\left[\rvcx_t\rvcx_t^T\right]\right]\nonumber\\
&\ge 0,\label{eq:stochastic_H_lower_bound}
\end{align}
where $\bar{\rvu}_t^{i}$ is the stochastic optimal control input signal of the agent $i$ at time instant $t$.
Next, by using \eqref{eq:stochastic_H_lower_bound} and following along the lines of the proof of Lemma \ref{lem:H_lower_bound}, it can be seen that $\left(\bar{Q},\{\bar{P}_t\}\right)$ is an optimal solution to \eqref{eq:opt_noise}, and that
\begin{align}
\mE \left[H_S^{\rvcy}(\bar{Q},\{\bar{P}_t\})\right] + \frac{1}{M} \sum_{i = 1}^M \sum_{t=1}^{N-1}\mE\left[\| \bar{\rvu}_t^{i}\|^2\right]=0.\label{eq:lower_bound_attained_noisy}
\end{align}
Therefore, what remains is to show the uniqueness of the optimal solution. To this end, by Assumption \ref{ass:persistent_excitation_noisy} the covariance matrix $\cov(\rvx_1^{i},\rvx_1^{i})$ is strictly positive definite, and therefore the second-order moment of the initial value $\mE[\rvx_1^{i} \rvx_1^{iT}] = \cov(\rvx_1^{i}, \rvx_1^{i})+ \mE[\rvx_1^{i}]\mE[\rvx_1^{i}]^T\succ 0$. By the fact that $A_{cl}(t;Q)$ is invertible for all $Q\in\mathbb{S}^n_+$, it holds that
\begin{align*}
&\mE[\rvx_t^{i} \rvx_t^{iT}] = \left[\prod_{k=1}^{t-1} A_{cl}(k;Q)\right]\mE[\rvx_1^{i} \rvx_1^{iT}]\left[\prod_{k=t-1}^{1} A_{cl}^T(k;Q)\right]\succ 0,
\end{align*}
cf.~Lemma~\ref{lem:persistent_excitation}, and hence $\mE[\rvcx_t\rvcx_t^T] \succ 0$.
Now suppose that there exists some other $(Q^*,\{P^*_t\})$ that is also optimal to \eqref{eq:opt_noise}. By \eqref{eq:stochastic_H_lower_bound}, \eqref{eq:lower_bound_attained_noisy},  we have
\begin{align*}
&0 = \mE[H_S^{\rvcy}(Q^*,\{P^*_t\})] + \frac{1}{M} \sum_{i = 1}^M \sum_{t=1}^{N-1}\mE\left[\| \bar{\rvu}^i_t \|^2\right]\\
&\ge \frac{1}{M}\sum_{t=1}^{N-1} \trace \big[(K_t^*-\bar{K}_t)^T(B^TP_{t+1}B+I)\\
&\quad \times (K_t^*-\bar{K}_t)\mE\left[\rvcx_t\rvcx_t^T\right]\big]\\
&=\frac{1}{M}\sum_{t=1}^{N-1}\left\|(B^TP_{t+1}B+I)^{\frac{1}{2}}(K_t^*-\bar{K}_t)\mE\left[\rvcx_t\rvcx_t^T\right]^{\frac{1}{2}}\right\|_F^2\\
&\ge 0,
\end{align*}
which together with the fact that $\mE[\rvcx_t\rvcx_t^T]\succ 0$ implies that
\begin{align*}
(B^TP_{t+1}B+I)^{\frac{1}{2}}(K_t^*-\bar{K}_t)\mE\left[\rvcx_t\rvcx_t^T\right]^{\frac{1}{2}}=0,
\end{align*}
for $t=1,\ldots,N-1$. The latter in turn implies that
\begin{align*}
K_t^*-\bar{K}_t=0,\forall t=1,\ldots,N-1,
\end{align*}
where $K_t^*$ denote the same expression as the one in the proof of Theorem \ref{thm:main_result_no_noise}.
To this end, following the same analysis as in the proof of Theorem \ref{thm:main_result_no_noise}, we can conclude that $Q^*=\bar{Q}$ and hence the statement holds.
\end{pf}

Since the distribution of the initial values $\rvx_1^{i}$ and the additive noise $\rvv_t^{i}$ are not known, we can not express the expected value of the objective function \eqref{eq:def_H_stochastic} explicitly and hence it is not possible to solve \eqref{eq:opt_noise} directly. Therefore, we derive an empirically estimate of the expectation based on the observations. To this end, first recall that due to Assumption~\ref{ass:noise_iid_and_zero_mean}, $\rvy_t^i$ are I.I.D for $i=1,\ldots,M$. This means that we can rewrite the cost \eqref{eq_opt_noise_cost} as 
\begin{align}
& \mE\left[ H_S^{\rvcy}(Q,\{P_t\}) \right] = \nonumber \\
&\mE\Big[ -\trace(P_{1} \rvy_1 \rvy_1^T) + \trace( P_N \rvy_N \rvy_N^{T}) +\sum_{t=1}^{N-1} \trace (Q \rvy_t \rvy_t^{T})\Big] \nonumber \\
& + \trace \left( P_1\Sigma \right) - \trace(P_N \Sigma) - (N-1)\trace (Q\Sigma), \label{eq:H_stochastic_rewritten}
\end{align}
where, for $t = 1, \ldots, N$, $\rvy_t$ is a random variable with the same distribution as $\rvy_t^{i}$ for $i = 1, \ldots, M$. 
This means that an empirical estimate of the expectation of \eqref{eq:def_H_stochastic} can be obtained as
\begin{align*}
& \mathbb{E} \left[\left( -\rvy_1^{T}P_1 \rvy_1 + \rvy_N^{T} P_N \rvy_N + \sum_{t= 1}^{N-1} \rvy_t^{T} Q \rvy_t \right) \right] \\
& \approx \frac{1}{M} \sum_{i = 1}^M \left( -\rvy_1^{iT}P_1 \rvy_1^{i} + \rvy_N^{iT} P_N \rvy_N^{i} + \sum_{t= 1}^{N-1} \rvy_t^{iT} Q \rvy_t^i \right).
\end{align*}

Based on this, we formulate the estimation problem
\begin{subequations}\label{eq:opt_noise_empirical}
  \begin{align}
    & \underset{Q,  \{P_t \}_{t=1}^N }{\text{minimize}}
    & &  H_{E}^{\rvcy}(Q,\{P_t\}) \\
    & \text{subject to}
    & & (Q,  \{P_t \}_{t=1}^N) \in \mathscr{D}(\varphi),
  \end{align}
where 
\begin{align}\label{eq:def_H_stochastic_empirical}
& H_{E}^{\rvcy}(Q,\{P_t\}) := \frac{1}{M}\Big[-\trace(P_{1}\rvcy_1^\phi \rvcy_1^{\phi T}) \nonumber \\
& \; + \trace(P_N \rvcy_N^\phi \rvcy_N^{\phi T}) +\sum_{t=1}^{N-1} \trace (Q \rvcy_t^\phi \rvcy_t^{\phi T}) \Big]\nonumber \\
& \; + \big( \trace \left( P_1\Sigma \right) - \trace(P_N \Sigma) - (N-1)\trace (Q\Sigma) \big).
\end{align}
\end{subequations}
The problem \eqref{eq:opt_noise_empirical} defines the estimator, and for a given realization $\{ Y^\phi_t \}_{t=1}^N$ of the stochastic variables $\{ \rvcy_t^{\phi} \}_{t = 1}^N$, it can be solved in order to obtain an estimate. In particular, we use the notation $H_{E}^{\rvcy}(Q,\{P_t\})|_{\rvcy=Y}$ to denote the cost function \eqref{eq:def_H_stochastic_empirical} evaluated at a particular realization. We now want to show that this estimator is in fact (asymptotically) statistically consistent.
However, note that since we approximate the expected value in the objective function by the empirical average, the objective function changes and hence the ``bounded-from-below" argument \eqref{eq:stochastic_H_lower_bound} does not necessarily hold for $H_{E}^{\rvcy}(Q,\{P_t\})|_{\rvcy = Y}$ on the domain $\mathscr{D}(\varphi)$.
This issue needs to be addressed in order to make \eqref{eq:opt_noise_empirical} well-posed. This is an important first step towards showing that the estimator is statistically consistent.

\begin{lemma}\label{lem:approx_problem_bounded_from_below}
The domain $\mathscr{D}(\varphi)$ in \eqref{eq:domain_D_varphi} is compact, and $H_{E}^{\rvcy}(Q,\{P_t\})|_{\rvcy=Y}$ is bounded on $\mathscr{D}(\varphi)$.
\end{lemma}
\begin{pf}
Consider the domain $\mathscr{D}(\varphi)$ and recall that, by the property of Schur complement, \eqref{eq:schur_complement_LMI} holds on the feasible domain \eqref{eq:LMI} and \eqref{eq:P_N_constr}.  Since the Frobenius norm is monotone with respect to the Loewner partial order,  it holds that
\begin{align*}
&\|P_t\|_F\le \\
&\|A^TP_{t+1}A+Q-A^TP_{t+1}B(B^TP_{t+1}B+I)^{-1}B^TP_{t+1}A\|_F
\end{align*}
By the Cauchy-Schwarz and the triangular inequality, we have that
\begin{align}
&\|P_t\|_F\le \|A\|_F^2\cdot \|P_{t+1}\|_F+\|Q\|_F \label{eq:norm_ineq_riccati} \\
&\qquad +\|A\|_F^2\cdot\|B\|_F^2\cdot\|P_{t+1}\|_F^2\cdot\|(B^TP_{t+1}B+I)^{-1}\|_F. \nonumber
\end{align}
Next, since $P_{t} \in \mS^n_+$, it holds that $B^TP_tB+I \succeq I$, and hence that $ (B^TP_tB+I)^{-1} \preceq I$ \citep[Cor.~7.7.4]{horn2013matrix}. By monotonicity of  the Frobenius norm with respect to the Loewner partial order, it therefore follows that $\|(B^TP_tB+I)^{-1}\|_F \le \|I\|_F = \sqrt{n}$. Now,  since $\|Q\|_F^2\le \varphi$ and $P_{N} = Q$, using this together with \eqref{eq:norm_ineq_riccati}, it follows that $\|P_{N-1}\|_F$ is bounded. Recursively applying this backwards for the time indices $t$, it follows that $\|P_t\|_F$ is bounded for all $t=1,\ldots,N$. This implies that the domain $\mathscr{D}(\varphi)$ is compact. Finally, since $H_{E}^{\rvcy}(Q,\{P_t\})|_{\rvcy=Y}$ is continuous, it is bounded on $\mathscr{D}$.
\end{pf}
\begin{remark}
Note that Assumption \ref{ass:Q_set_compact} is critical in the proof of Lemma \ref{lem:approx_problem_bounded_from_below}, since we can thus optimize over
$\mathscr{D}(\varphi)$ instead of $\mathscr{D}$. In fact,  $H_{E}^{\rvcy}(Q,\{P_t\})|_{\rvcy = Y}$ might not be bounded from below if we only impose $Q\in\mS^n_+$. To see this, assume that there exists a realization $\{Y_t^{\phi}\}$ of $\{\rvcy_t^{\phi}\}$ such that $Y_1^\phi Y_1^{\phi T} - M\Sigma=0$ and $Y_t^\phi Y_t^{\phi T}-M\Sigma\prec 0$, $t=2,\ldots,N$. Let $(Q,\{P_t\})$ satisfies DRE \eqref{eq:DRE}, and note that then $(\alpha Q,\{\alpha P_t\})$ also satisfies DRE \eqref{eq:DRE} for any positive $\alpha$.  Therefore, $(\alpha Q,\{\alpha P_t\}) \in \mathscr{D}$ for all $\alpha >0$, and thus $\mathscr{D}$ is not a bounded set. Moreover, as $\alpha\rightarrow \infty$ it holds that
\begin{align*}
&H_{E}^{\rvcy}(\alpha Q,\{\alpha P_t\})|_{\rvcy=Y} = \frac{1}{M} \Bigg[ - \trace\left(\alpha P_1(Y_1^\phi Y_1^{\phi T}-M\Sigma)\right)\\
& + \trace\left(\alpha P_N(Y_N^\phi Y_N^{\phi T}-M\Sigma)\right)\\
& + \sum_{t=1}^{N-1} \trace\left(\alpha Q(Y_t^\phi Y_t^{\phi T}-M\Sigma)\right) \Bigg] \rightarrow -\infty,
\end{align*}
hence $H_{E}^{\rvcy}(\alpha Q,\{\alpha P_t\})|_{\rvcy=Y}$ is not bounded from below on $\mathscr{D}$.
\end{remark}

Next, we show that the ``Uniform Law of Large Numbers" holds for $H_{E}^{\rvcy}(Q,\{P_t\})$. 

\begin{lemma}[Uniform law of large numbers]\label{lem:uniform_law_of_large_number}
Under Assumption \ref{ass:noise_iid_and_zero_mean}, \ref{ass:persistent_excitation_noisy} and \ref{ass:Q_set_compact}, it holds that
  \begin{align*}
  \sup_{(Q,\{P_t\}_{t=1}^N) \in \mathscr{D}(\varphi)} \big|H_{E}^{\rvcy}(Q,\{P_t\}) - \mE\left[H_S^{\rvcy}(Q,\{P_t\})\right]\big|\overset{a.s.}{\rightarrow} 0,
  \end{align*}
as $M \to \infty$.
\end{lemma}
\begin{pf}
It is clear that $H_S^{\rvcy}(Q,\{P_t\})$ is continuous with respect to $\{ Y^\phi_t \}_{t=1}^N$ and therefore it is a measurable function of $\{ Y^\phi_t \}_{t=1}^N$ for each $Q$ and $\{P_t\}$. 
On the other hand, Assumption \ref{ass:persistent_excitation_noisy} implies that $\mE[\|\rvx_t^{i} \|^2]<\infty$ for $i = 1,\ldots, M$ \citep[cf.~the proof of Theorem 4.1]{zhang2019inverse}. By Assumption \ref{ass:noise_iid_and_zero_mean}, we have $\mE[\|\rvv_t^{i} \|^2]<\infty$.  Since $\rvx_t^{i}$ is independent of $\rvv_t^{i}$, and since $\mE[\rvv_t^{i}]=0$, it follows that
\begin{align*}
\mE[\|\rvy_t^{i}\|^2]&=\mE[(\rvx_t^{i}+\rvv_t^{i})^T(\rvx_t^{i}+\rvv_t^{i})]=\mE[\|\rvx_t^{i}\|^2]+\mE[\|\rvv_t^{i}\|^2]\\
&< \infty.
\end{align*}
In addition, from Lemma \ref{lem:approx_problem_bounded_from_below}, we know that there exists constants $\{\bar{\varphi}_t\}_{t=1}^N$ such that for all $(Q,\{P_t\}_{t=1}^N)\in\mathscr{D}(\varphi)$ we have that $\|P_t\|_F\le \bar{\varphi}_t$. 
Using the form of \eqref{eq_opt_noise_cost} given in \eqref{eq:H_stochastic_rewritten}, by the Cauchy-Schwarz and the triangular inequality it therefore holds that
\begin{align*}
&H_S^{\rvcy}(Q,\{P_t\}) \\
&= - \trace(P_1\rvy_1\rvy_1^T) + \trace(P_N\rvy_N\rvy_N^T)+\sum_{t=1}^{N-1}\trace(Q\rvy_t\rvy_t^T)\\
&\quad + \trace(P_1\Sigma)-\trace(P_N\Sigma)-(N-1)\trace(Q\Sigma)\\
&\le \left|\rvy_1^TP_1\rvy_1\right|+\left|\rvy_N^TP_N\rvy_N\right|+\sum_{t=1}^{N-1}\left|\rvy_t^TQ\rvy_t\right|\\
&\quad +\left| \trace(P_1\Sigma)\right|+\left| \trace(P_N\Sigma)\right|+(N-1)\left|\trace(Q\Sigma)\right|\\
&\le \|\rvy_1\|^2\|P_1\|_F+\|\rvy_N\|^2\|P_N\|_F+\sum_{t=1}^{N-1}\|\rvy_t\|^2\|Q\|_F\\
&\quad +\|P_1\|_F\|\Sigma\|_F+\|P_N\|_F\|\Sigma\|_F+(N-1)\|Q\|_F\|\Sigma\|_F\\
&\le \bar{\varphi}_1(\|\rvy_1\|^2+\|\Sigma\|_F)+\bar{\varphi}_N(\|\rvy_N\|^2+|\Sigma\|_F)\\
&\quad+\varphi\sum_{t=1}^{N-1}(\|\rvy_t\|^2+\|\Sigma\|_F):=d(\{\rvy_t\}),
\end{align*}
and it is clear that $\mE[d(\{\rvy_t\})]<\infty$ since $\mE[\|\rvy_t^{i} \|^2]<\infty$. Therefore, by \citep[Thm.~2]{jennrich1969asymptotic} the result follows.
\end{pf}

We are now ready to prove the main result of this section.

\begin{theorem}[Statistical consistency]\label{thm:statistical_consistency}
Suppose that $(Q_M^*,\{P_{t,M}^*\}_{t=1}^N)$ is an optimal solution to \eqref{eq:opt_noise_empirical} when observing $M$ agents. Then $Q_M^*\overset{p}\rightarrow \bar{Q}$ as $M\rightarrow \infty$, where $\bar{Q}$ is the true parameter used in the objective function of ``forward" problem \eqref{eq:stochastic_forward_problem}.
\end{theorem}

\begin{pf}
The theorem is proved by showing that all the conditions in \cite[Thm.~5.7]{van1998asymptotic} are satisfied.
To this end, the first condition follows from Lemma \ref{lem:uniform_law_of_large_number}, since convergence a.s.~implies convergence in probability \citep[Lem~3.2]{kallenberg1997foundations}.
Next, the second condition holds since by Proposition~\ref{prop:stoc_prob_unique_sol} the optimal solution to \eqref{eq:opt_noise} is unique, together with the fact that $\mathcal{D}$ is compact (see \cite[p.~46]{van1998asymptotic}).
Therefore, all conditions in \cite[Thm.~5.7]{van1998asymptotic} are satisfied, and the statement hence follows.
\end{pf}

\section{On numerical ill-conditioning}\label{sec:ill-cond}
Proposition~\ref{prop:global_identifiability} shows that the model is globally identifiable, and Theorem~\ref{thm:main_result_no_noise} shows that the optimization problem \eqref{eq:opt_no_noise} has a unique optimal solution at the ``true" $\bar{Q}$. Hence, in theory the latter can be recovered by solving the optimization problem. 
Nevertheless, recovering this optimal solution turns out to be numerically difficult, in particular for certain problem instances - this will be demonstrated with examples in Section~\ref{sec:numerics}. Here, we argue that this has to do with an intrinsic numerical ill-conditioning of the inverse problem for these problem instances.

To this end, recall that Pontryagin's Maximum Principle (PMP) gives a necessary and sufficient condition for optimality in the forward problem \eqref{eq:forward_problem}.  Namely, $\{x_t\}_{t=1}^N$ and $\{u_t\}_{t=1}^{N-1}$ are the optimal trajectory and control signal of \eqref{eq:forward_problem}, respectively, if and only if there exists adjoint states $\{\lambda_t\}_{t=2}^N$ such that
\begin{subequations}\label{eq:PMP_original}
\begin{align}
&\lambda_t=A^T\lambda_{t+1}+Qx_t,\:t=2,\ldots,N-1,\label{eq:PMP_original_1}\\
&\lambda_N=Qx_N,\label{eq:PMP_original_2}\\
&u_t=-B^T\lambda_{t+1},\:t=1,\ldots,N-1. \label{eq:PMP_original_3}
\end{align}
\end{subequations}
Based on \eqref{eq:PMP_original_1}--\eqref{eq:PMP_original_2}, we can write a linear system of equations for the adjoint variables of the $i$th agent, namely
\[
\begin{bmatrix}
I &-A^T\\
  &I &\ddots\\
  &  &\ddots & -A^T\\
  & & &I
\end{bmatrix}
\underbrace{\begin{bmatrix}
\lambda_2^i\\\vdots\\\lambda_N^i
\end{bmatrix}}_{=: \lambda^i}=(I\otimes Q)
\underbrace{\begin{bmatrix}
x_2^i\\\vdots\\x_N^i
\end{bmatrix}}_{=: x_{2:N}^i}.
\]
Solving this for the adjoint variables gives
\[
\lambda^i
=
\begin{bmatrix}
I &A^T &(A^T)^2 &\cdots &(A^T)^{N-2}\\
  &I &A^T&\cdots &(A^T)^{N-3}\\
  &  &\ddots &\ddots &\vdots\\
  &  &      &I&A^T\\
  &  &      & &I
\end{bmatrix}(I\otimes Q)
x_{2:N}^i,
\]
and substituting the latter into \eqref{eq:PMP_original_3} we obtain
\begin{align}
&-\underbrace{\begin{bmatrix}
u_1^i\\\vdots \\u_{N-1}^i
\end{bmatrix}}_{=: u^i}
=(I_{N-1} \otimes B^T)\lambda^i\nonumber\\
&= \! \underbrace{\begin{bmatrix}
B^T \! &B^TA^T \! &B^T(A^T)^2 \! &\cdots \! &B^T(A^T)^{N-2}\\
      &B^T \!     &B^TA^T \!     &\cdots \! &B^T(A^T)^{N-3}\\
      &           &\ddots \!     &\ddots \! &\vdots\\
      &           &              &B^T \!    &B^TA^T\\
      &           &              &          &B^T
\end{bmatrix}}_{= [\mathcal{S}_2\Gamma,\cdots,\mathcal{S}_N\Gamma]=\mathcal{S}(I_{N-1}\otimes\Gamma)} \! (I_{N-1} \otimes Q)
x_{2:N}^i\label{eq:PMP_u},
\end{align}
where $\mathcal{S} := [\mathcal{S}_2,\cdots,\mathcal{S}_N]$, $\mathcal{S}_N := I_{m(N-1)}$, $\mathcal{S}_{N-k} \in \mR^{m(N-1) \times m(N-1)}$ is a block-matrix with identity matrices of size $n \times n$ on the $k$th upper block-diagonal, and
\[
\Gamma :=
\begin{bmatrix}
(A)^{N-2}B & (A)^{N-3}B & \cdots & B
\end{bmatrix}^T.
\]
Using \eqref{eq:PMP_u} and the fact that
\begin{align*}
x_t^i = \left[A^{t-2}B\: \cdots\:AB\: B\right] \left[u_1^{iT},\cdots,u_{t-1}^{iT}\right]^T+A^{t-1}x_1^i,
\end{align*}
we have
\begin{align*}
&x_{2:N}^i = \underbrace{\begin{bmatrix}
B&0&\cdots&0&0\\
AB&B&\cdots&0&0\\
\vdots&\vdots & & \vdots &\vdots\\
A^{N-1}B&A^{N-2}B&\cdots&AB&B
\end{bmatrix}}_{= (I_{N-1}\otimes \Gamma^T)\mathcal{S}^T}u^i+
\underbrace{\begin{bmatrix}
A\\A^2\\\vdots\\A^{N-1}
\end{bmatrix}}_{=: \tilde{A}}x_1^i,\nonumber\\
&=-\underbrace{(I_{N-1}\otimes \Gamma^T)\mathcal{S}^T\mathcal{S}(I_{N-1}\otimes \Gamma)(I_{N-1}\otimes Q)}_{=: \mathscr{F}(Q)}x_{2:N}^i+\tilde{A}x_1^i\\
&\implies(I_{(N-1)n}+\mathscr{F}(Q))x_{2:N}^i = \tilde{A}x_1^i
\end{align*}
Since \eqref{eq:forward_problem} has a unique solution, $I_{(N-1)n}+\mathscr{F}(Q)$ is intrinsically invertible. Therefore, $(I_{(N-1)n}+\mathscr{F}(Q))^{-1}\tilde{A}$ is another representation of the same model structure as $\mathcal{M}(Q, \{ \mathcal{P}_t \}_{t = 2}^N)$, if the former is also combined with the permutation matrices $\{\mathscr{P}_t\}_{t=2}^N$; for the sake of brevity we omit the details.

From above, it can be seen that the model structure is identifiable at $\bar{Q}$ if and only if
\begin{align}
\Upsilon \! := \! \left\{\img(\Delta Q)\mid \Delta Q\in\mathbb{S}^n,\bar{Q}+\Delta Q\in\mathbb{S}^n_+\right\}\cap \kernl(\Gamma) \! = \! \{0\}.
\label{eq:model_identifiability_condition}
\end{align}
Assuming that the system $(A,B)$ is controllable implies that the controllability matrix
\begin{align*}
\Gamma_n = \begin{bmatrix}
(A)^{n-1}B &\cdots &AB &B
\end{bmatrix}^T
\end{align*}
has full column-rank, and hence that $\kernl(\Gamma)=\{0\}$. The latter, in turn, means that \eqref{eq:model_identifiability_condition} is fulfilled. Nevertheless, in practice, if $\Gamma$ is ill-conditioned, the kernel of $\Gamma$ can be ``expanded" from a numerical perspective. In this case we can have $\Gamma \Delta Q \approx 0$ for some $\Delta Q$ that is not close to zero, and the set $\Upsilon$ might numerically not be the singleton $\{ 0 \}$. Thus, it is possible that
\[
(I_{(N-1)n}+\mathscr{F}(\bar{Q} + \Delta Q))^{-1}\tilde{A} \approx (I_{(N-1)n}+\mathscr{F}(\bar{Q}))^{-1}\tilde{A}
\]
for some feasible $\Delta Q$ that is not close to zero. Therefore, when the controllability matrix is ill-conditioned, we might have very similar models that corresponds to very different $Q$'s. In general, we expect that it will be numerically challenging to recover the ``true'' $\bar{Q}$ in these settings, regardless of which method is used.
However, although it might be numerically difficult to recover $\bar{Q}$ in these circumstances, the proof of Lemma~\ref{lem:H_lower_bound} ensures that the control gains corresponding to $\bar{Q} + \Delta Q$ numerically coincides with the ``true" control gain, which is sufficient for predicting the agents' behaviors. Finally, a similar argument holds in the case of noisy observations.

\section{Numerical experiments and discussions}\label{sec:numerics}
In this section, we present a number of numerical experiments, performed on a number of different discrete-time systems, to illustrate the properties of the proposed algorithm.
In particular, the discrete-time systems are all generated by sampling continuous-time systems $\dot{x} = \hat{A}x + \hat{B}u$ via $A = e^{\hat{A}\Delta t}$ and $B =  \int_0^{\Delta t} e^{\hat{A}t} dt \hat{B}$, where the sampling period $\Delta t=0.05$.

All numerical examples are run on a MacBook Pro with Apple M1 eight-core CPU and 16GB of RAM.
The solutions are obtained by implementing the optimization problems in Matlab using YALMIP \citep{lofberg2004yalmip} and solving them using MOSEK \citep{mosek}.

\subsection{Noiseless case}\label{subsec:num_noiseless}
\begin{figure}[!htpb]
    \centering
    \includegraphics[width=0.9\columnwidth]{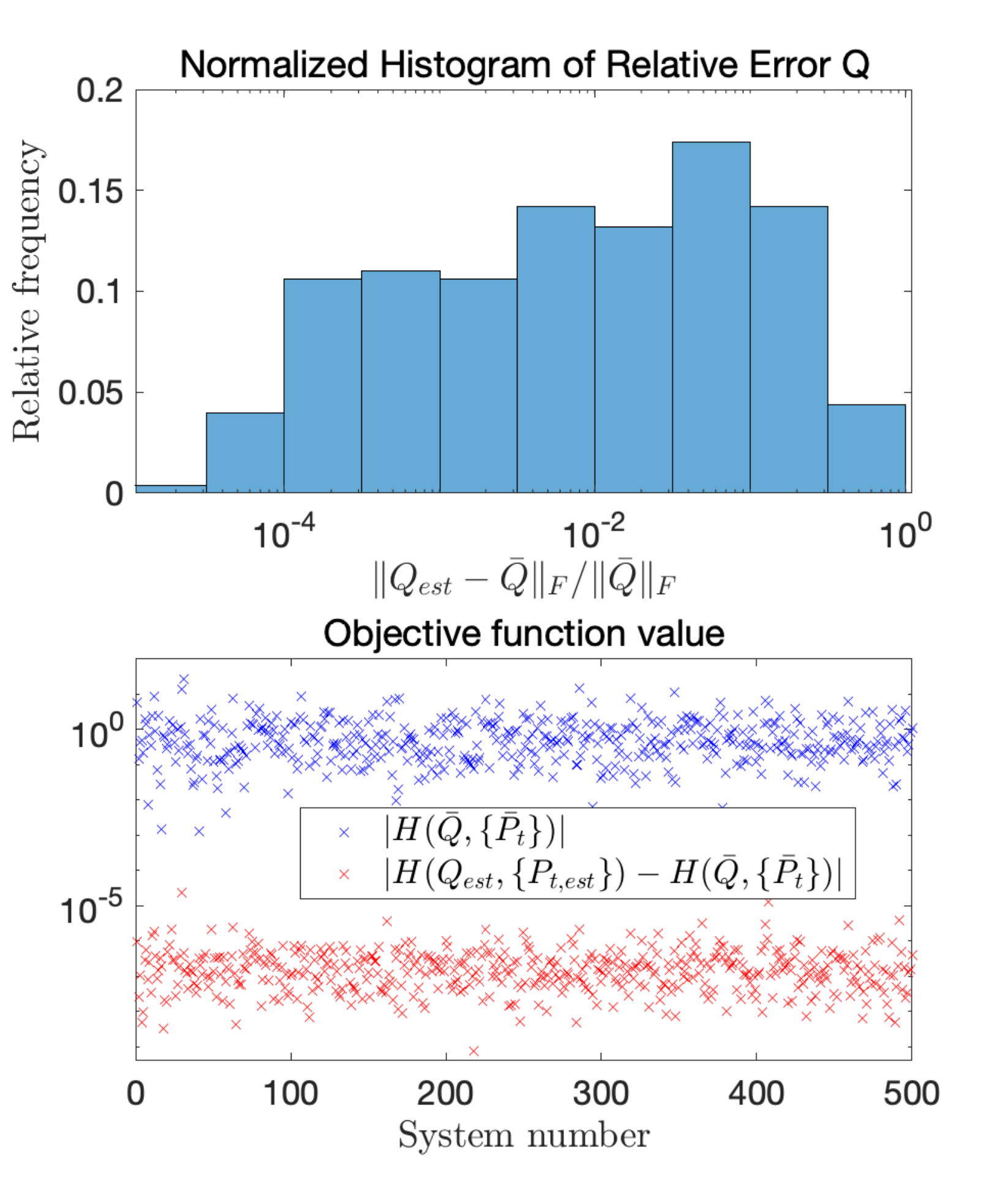}
    \caption{The upper plots shows a (normalized) histogram of the relative error in the estimate $Q_{est}$ obtained with noiseless data, as described in Section~\ref{subsec:num_noiseless}. The lower plots shows the absolute value of the (scaled) gobally optimal objective function value, as well as absolute value of the difference between the globally optimal objective function value and the objective function value in the obtained estimate.}
    \label{fig:noiseless-1}
\end{figure}

\begin{figure}[!htpb]
    \centering
    \includegraphics[width=0.9\columnwidth]{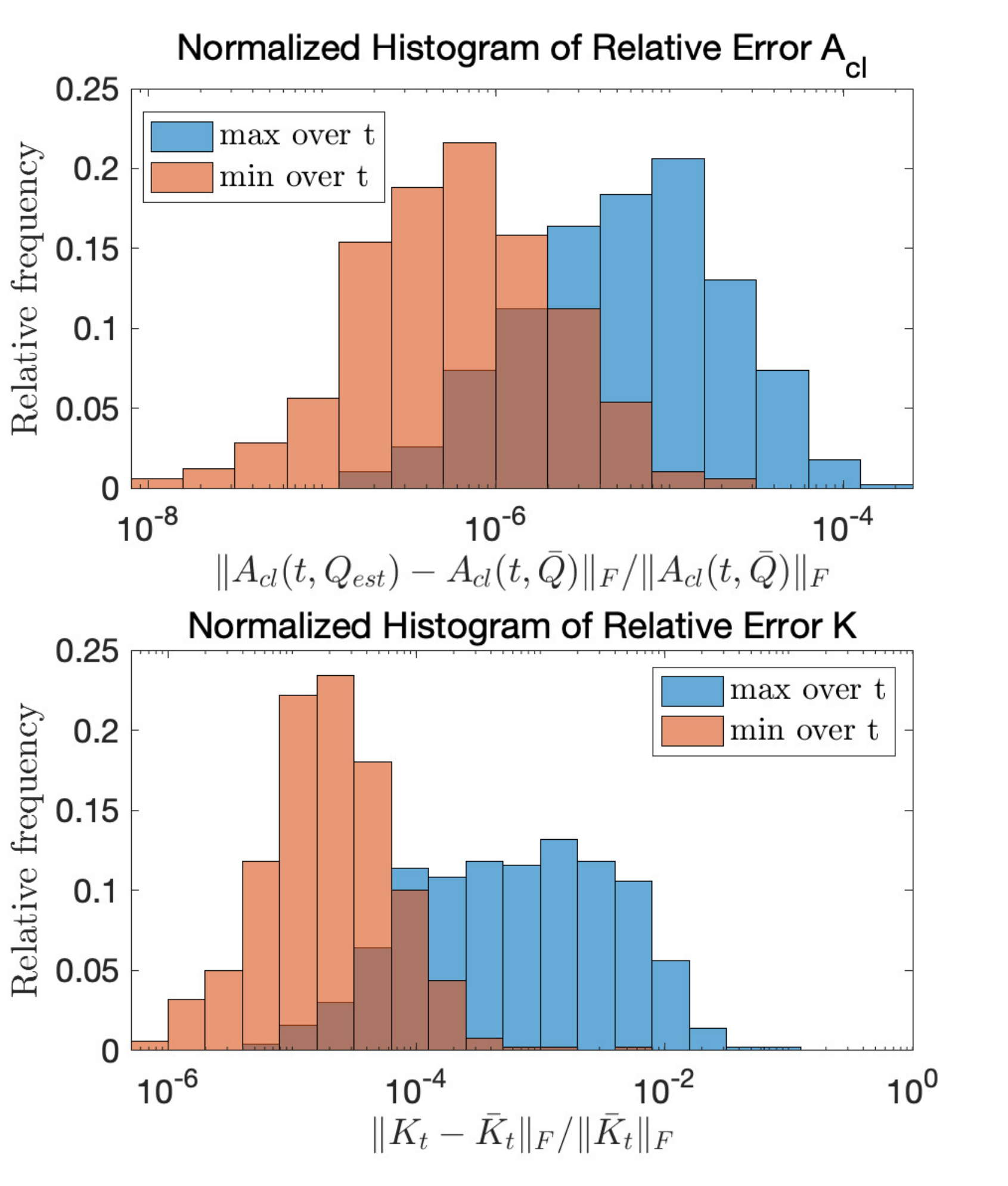}
    \caption{Histograms (normalized) of the relative error in the estimates of the closed-loop system matrix $A_{cl}(t, Q_{est})$ and the gain $K_t$, obtained from noiseless data as described in Section~\ref{subsec:num_noiseless}. Since the quantities are time-varying, the largest and smallest relative errors are shown.}
    \label{fig:noiseless-2}
\end{figure}

In this experiment, the dimension of the system is set to $n = 3$, and $m = 1$. More specifically, we randomly generate system matrices $\hat{A} \in \mR^{3 \times 3}$ and $\hat{B} \in \mR^{3 \times 1}$ with entries drawn from a normal distribution with mean value zero and standard deviation one, i.e., with entries drawn from the distribution $\mathcal{N}(0,1)$. These are then sampled to generate discrete-time systems, as described above. Moreover,
the ``true "$\bar{Q}$ is randomly generated as $\bar{Q} = GG^T$, where $G \in \mR^{3 \times 3}$ with entries drawn from $\mathcal{N}(0,1)$. We let $\bar{Q} \in \bar{\mathbb{S}}^3_+(5)$; any randomly generated $\bar{Q}$ would be discarded if it does not belong to $\bar{\mS}^3_+(5)$ and another random $\bar{Q}$ would be generated.
In this way, 500 random triplets $(A, B, \bar{Q})$ are generated.
For each such random triplet, we set the time horizon to $N = 20$ and generate $M = 15$ random starting points $x_1^i$.
The latter are drawn from a uniform distribution on $[-10, 10] \times [-10, 10] \times [-10, 10]$.
The forward problem \eqref{eq:forward_problem} is solved for each starting point, and the noiseless data is then used to solve the optimization problem in \eqref{eq:opt_no_noise}, except that the cost function is scaled with $10^{-4}$ in order to give a better numerical scaling for the problem.
While this does not change any analytic properties, the obtained optimal solutions were observed to have a smaller relative error in general.
The results are presented in \figref{fig:noiseless-1} and \figref{fig:noiseless-2}.

The lower plot in \figref{fig:noiseless-1} illustrates the absolute value of the scaled objective function value at the theoretical optimal solution $(\bar{Q}, \{\bar{P}_t\})$, together with the absolute value of the difference between the objective function value at $(\bar{Q}, \{\bar{P}_t\})$ and at the solution $(Q_{est}, \{ P_{t, est} \})$ obtained with the solver. 
As can be seen, the difference is in general several orders of magnitude smaller than the optimal value of the cost function, despite the fact that the obtained $Q_{est}$ is sometimes relatively far from $\bar{Q}$ (cf.~upper plot in \figref{fig:noiseless-1}).
This indicates that the cost function is ``flat" in a region around the optimal solution, which makes the ``true" $\bar{Q}$ hard to recover numerically with high accuracy. We believe that this is highly related to the fact that the controllability matrix is very ill-conditioned (the condition numbers of which varies from $2.899\times 10^2$ to $7.439\times 10^5$). Mitigation of this numerical difficulties is left for future work.

Nevertheless, as can be seen from \figref{fig:noiseless-2}, the corresponding control gain and closed loop system matrix are well-recovered, which serves the purpose of ``predicting the agent's behaviour". In fact, the latter is in general recovered with better accuracy than the former.
This indicates that for certain systems, a larger mismatch in $Q$ can still give small mismatches in the control gains and the closed-loop system matrices. Hence it is harder to identify $Q$ numerically in these cases.
Since the control gain and the closed-loop system are time-varying, the smallest and largest relative error over all time points are shown in \figref{fig:noiseless-2}.

The fact that the closed-loop system matrix in general seems to be better recovered than $\bar{Q}$ seems to indicate that the ``flatness'' of the cost function for certain problem instances is (at least partly) related to the discussion in Section~\ref{sec:ill-cond}. Namely, that for certain problem instances, substantially different $Q$'s can give rise to very similar closed-loop system matrices. 

\subsection{Noisy case}\label{subsec:num_noisy}

Next, we illustrate the statistical consistency of the method. This is done on a dynamical system which does not have a numerically ill-conditioned controllability Gramian. More specifically, the continuous-time dynamics is given by
\[
\hat{A} = \begin{bmatrix}
0 & 0 \\
0 & 0
\end{bmatrix},
\qquad
\hat{B} = \begin{bmatrix}
1 & 0 \\
0 & 1
\end{bmatrix},
\]
which is the kinematic dynamics of a point mass that moves on a two-dimensional plane. In fact, the condition number of the controllability Gramian for this dynamical system is $1$.
Moreover, with the corresponding discretized system matrices $(A,B)$, for any $\Qtrue \succ 0$ the forward problem \eqref{eq:forward_problem} describes an agent that moves towards the origin. Hence, a group of homogeneous agents that are all governed by \eqref{eq:forward_problem} is a model for the simplified setting of ``non-interacting'' agents moving towards a common goal at the origin. As mentioned in the introduction,  in the future we intend to extend the work to interacting agents.

\begin{figure}[!htpb]
    \centering
    \includegraphics[trim=1cm 1.7cm 1.5cm 0.5cm, clip, width = \columnwidth]{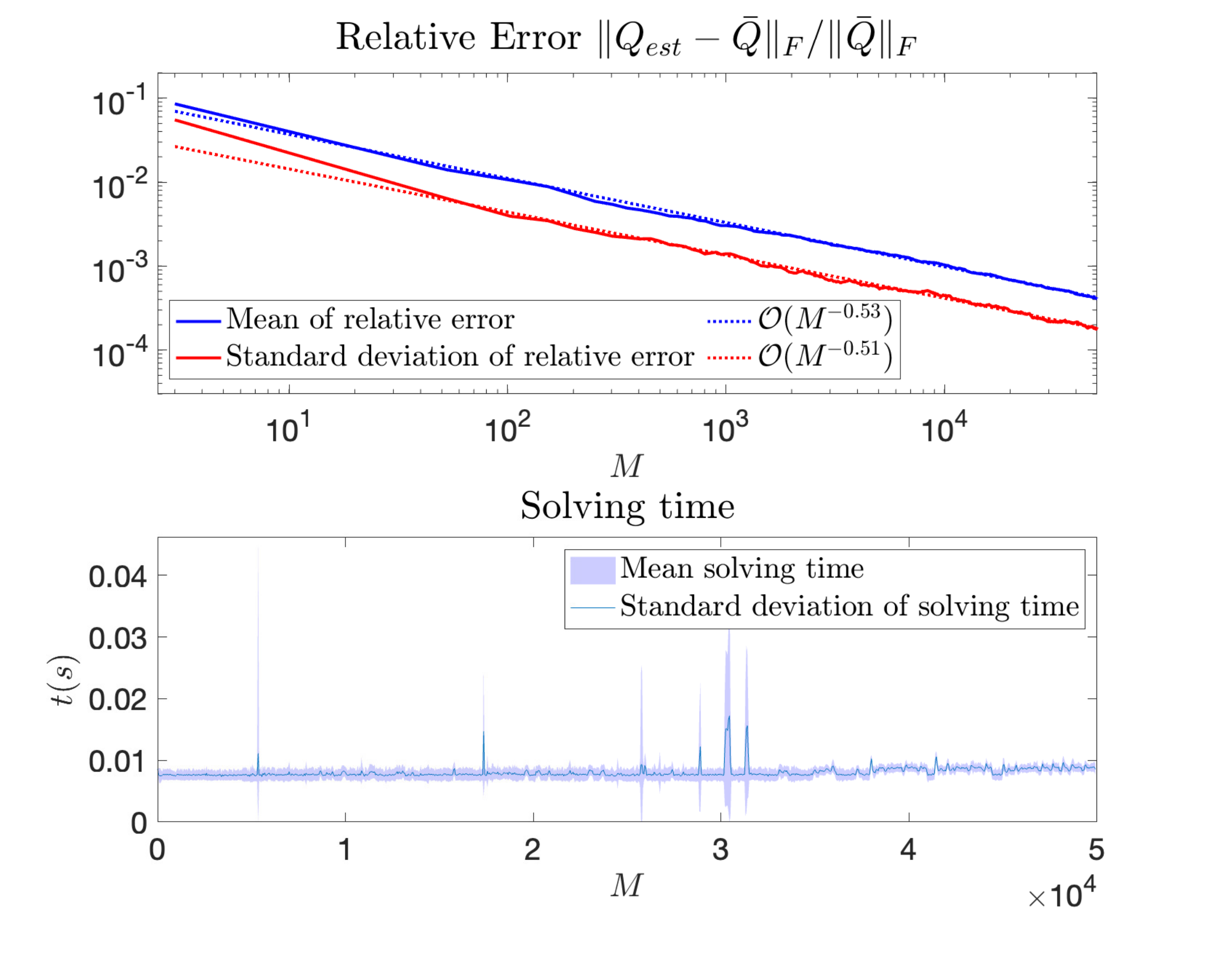}
    \caption{The upper plot shows the mean and standard deviation of the relative error of $Q_{est}$ as a function of the number of agents. The estimates are obtained using noisy data, as described in Section~\ref{subsec:num_noisy}. Moreover, the lower plot shows the time it took (in seconds) to solve the corresponding optimization problem.}
    \label{fig:noisy-2}
\end{figure}

The ``true'' $\Qtrue$ is generated as described in Section~\ref{subsec:num_noiseless}, and we set the time horizon to $N = 20$.
We then generate  $49 953$ random starting points $x_1^i$, drawn from a uniform distribution supported on $[-10, 10] \times [-10, 10]$.
For each agent, the forward problem \eqref{eq:forward_problem} is solved, and noise is added on the obtained optimal states (including the initial state).
The additive noise is drawn from a multi-variate zero-mean normal distribution with covariances matrix drawn from a Wishart distribution of degree $2$, i.e., with the same number of degrees of freedom as the dimension of the state space.
The Wishart distribution has a random covariance generated as $0.02 GG^T$, where each element in $G \in \mR^{2 \times 2}$ was drawn from $\mathcal{N}(0,1)$. 
The trajectories are then divided into groups of size $M = 3 +50 (k-1) $, for $k = 1, \ldots, 1000$,  where each larger group contains all the trajectories of a smaller group.
For the fixed triplet $(A, B, \bar{Q})$, the above process is repeated for $100$ times so that, in total, we get $100$ noisy data sets with a varying number of agents in each.
Moreover, the signal-to-noise ratio (SNR) in the data sets varies between $29.2479$ and $29.3767$ dB.%
\footnote{The SNR in a data set is computed as the mean of the SNR for all trajectories in that data set, where the SNR of a trajectory is computed as the sum of the squared norm of all states divided by the sum of the squared norm of all noise realizations.}
For each data set and each trajectory number $M$,  which varies from $3$ to $49953$, the problem \eqref{eq:opt_noise_empirical} is solved, using the corresponding cost function $H_S^{(\rvcy, M)}(Q,\{P_t\})|_{\rvcy=Y}$. 
That means that for each fixed $M$, we get one hundred estimates of $\bar{Q}$, and from these one hundred estimates we calculate the mean and the standard deviation of the relative error $\| Q_{est} - \Qtrue \|_F/ \| \Qtrue \|_F$. The result, as a function of $M$, are shown in \figref{fig:noisy-2}.

From the upper plot in \figref{fig:noisy-2} we see that both the mean and the standard deviation of the relative error of the estimates decreases with increasing $M$, in line with the statistical consistency of the estimate as proved in Theorem~\ref{thm:statistical_consistency}. Moreover, in the log-log plot of the mean and the standard deviation of the estimates v.s.~$M$, we can see that the relation is approximately linear. Fitting a log-linear model to the data, i.e., fitting an affine function to the logarithmic data, we get that $\texttt{Mean of relative error} \approx \mathcal{O}(M^{-0.53})$ and $\texttt{Standard deviation of relative error} \approx \mathcal{O}(M^{-0.51})$. The corresponding lines are also shown in  \figref{fig:noisy-2}. The orders are close to $-0.5$, and hence we suspect that the convergence rate is $\mathcal{O}(M^{-0.5})$ and that $\sqrt{M}(Q_M-\Qtrue)$ is asymptotically normal, just like most M-estimators such as maximum log-likelihood \citep[p.~51]{van1998asymptotic}. Further analysis of this is left for future work.

Finally, from the lower plot in \figref{fig:noisy-2} we see that, as expected from Remark~\ref{rem:scaling_num_agents}, the time to solve the corresponding optimization problem does note scale with the number of agents $M$.%
\footnote{The solving times are returned by YALMIP.}

\section{Conclusions}\label{sec:conclusion}
In this work we have considered the linear-quadratic inverse optimal control problem in discrete time and with finite time horizon, but where the observed homogeneous agents are indistinguishable.
In the case of exact measurements of the states, we show that the true parameter $\bar{Q}$ can be recovered as the unique globally optimal solution to a semidefinite programming problem.
Moreover, the size of this convex optimization problem is independent of the number of agents observed, and the formulation is thus suitable also for scenarios with a large number of agents.
Furthermore, in the case of noisy state observations the optimization problem is modified, and a statistically consistent estimator is obtained as the unique globally optimal solution to another semidefinite programming problem.
However, from numerical simulations it seems that that for certain problem instances the cost functions in both semidefinite programs are relatively flat around the globally optimal solutions, and hence accurate estimates of the parameters are difficult to recover.
Nevertheless, estimates of the time-varying control gains and closed-loop system matrices are obtained with higher accuracy.
An analysis of numerical ill-conditioning of the problem seem to suggest that this could be linked to the fact that for certain problem instances, different $Q$'s can give rise to similar close-loop system behaviour.
An important open question is therefore if and how the numerical conditioning of the estimators can be improved.

\bibliographystyle{plain}
\bibliography{} 

\begin{thebibliography}{10}

\bibitem{alexander1996optima}
R.~McNeill Alexander.
\newblock {\em Optima for animals}.
\newblock Princeton University Press, Princeton, NJ, 1996.

\bibitem{anderson2007optimal}
Brian D~O Anderson and John~B Moore.
\newblock {\em Optimal control: linear quadratic methods}.
\newblock Dover publications, Mineola, NY, 2007.

\bibitem{mosek}
MOSEK ApS.
\newblock {\em The MOSEK optimization toolbox for MATLAB manual. Version 9.0.},
  2019.

\bibitem{aswani2018inverse}
Anil Aswani, Zuo-Jun Shen, and Auyon Siddiq.
\newblock Inverse optimization with noisy data.
\newblock {\em Operations Research}, 66(3):870--892, 2018.

\bibitem{berret2011evidence}
Bastien Berret, Enrico Chiovetto, Francesco Nori, and Thierry Pozzo.
\newblock Evidence for composite cost functions in arm movement planning: an
  inverse optimal control approach.
\newblock {\em PLoS computational biology}, 7(10):e1002183, 2011.

\bibitem{boyd1994linear}
Stephen Boyd, Laurent El~Ghaoui, Eric Feron, and Venkataramanan Balakrishnan.
\newblock {\em Linear matrix inequalities in system and control theory}.
\newblock SIAM, Philadelphia, PA, 1994.

\bibitem{finn2016guided}
Chelsea Finn, Sergey Levine, and Pieter Abbeel.
\newblock Guided cost learning: Deep inverse optimal control via policy
  optimization.
\newblock In {\em International conference on machine learning}, pages 49--58.
  PMLR, 2016.

\bibitem{fridovich2020efficient}
David Fridovich-Keil, Ellis Ratner, Lasse Peters, Anca~D Dragan, and Claire~J
  Tomlin.
\newblock Efficient iterative linear-quadratic approximations for nonlinear
  multi-player general-sum differential games.
\newblock In {\em 2020 IEEE International Conference on Robotics and Automation
  (ICRA)}, pages 1475--1481. IEEE, 2020.

\bibitem{hatz2012estimating}
Kathrin Hatz, Johannes~P Schloder, and Hans~Georg Bock.
\newblock Estimating parameters in optimal control problems.
\newblock {\em SIAM Journal on Scientific Computing}, 34(3):A1707--A1728, 2012.

\bibitem{horn1994topics}
Roger~A. Horn and Charles~R. Johnson.
\newblock {\em Topics in matrix analysis}.
\newblock Cambridge University Press, New York, NY, 1994.

\bibitem{horn2013matrix}
Roger~A. Horn and Charles~R. Johnson.
\newblock {\em Matrix analysis}.
\newblock Cambridge university press, New York, NY, 2013.

\bibitem{jennrich1969asymptotic}
Robert~I Jennrich.
\newblock Asymptotic properties of non-linear least squares estimators.
\newblock {\em The Annals of Mathematical Statistics}, 40(2):633--643, 1969.

\bibitem{ji2019feedback}
Xuewu Ji, Kaiming Yang, Xiaoxiang Na, Chen Lv, Yulong Liu, and Yahui Liu.
\newblock Feedback game-based shared control scheme design for emergency
  collision avoidance: a fuzzy-linear quadratic regulator approach.
\newblock {\em Journal of Dynamic Systems, Measurement, and Control}, 141(8),
  2019.

\bibitem{jin2019inverse}
Wanxin Jin, Dana Kuli{\'c}, Jonathan Feng-Shun Lin, Shaoshuai Mou, and Sandra
  Hirche.
\newblock Inverse optimal control for multiphase cost functions.
\newblock {\em IEEE Transactions on Robotics}, 35(6):1387--1398, 2019.

\bibitem{kallenberg1997foundations}
Olav Kallenberg.
\newblock {\em Foundations of modern probability}.
\newblock Springer, 1997.

\bibitem{kalman1964linear}
Rudolf~E. Kalman.
\newblock When is a linear control system optimal?
\newblock {\em Journal of Basic Engineering}, 86(1):51--60, 1964.

\bibitem{keshavarz2011imputing}
Arezou Keshavarz, Yang Wang, and Stephen Boyd.
\newblock Imputing a convex objective function.
\newblock In {\em 2011 IEEE international symposium on intelligent control},
  pages 613--619. IEEE, 2011.

\bibitem{kopf2017inverse}
Florian K{\"o}pf, Jairo Inga, Simon Rothfu{\ss}, Michael Flad, and S{\"o}ren
  Hohmann.
\newblock Inverse reinforcement learning for identification in linear-quadratic
  dynamic games.
\newblock {\em IFAC-PapersOnLine}, 50(1):14902--14908, 2017.

\bibitem{li2020continuous}
Yibei Li, Yu~Yao, and Xiaoming Hu.
\newblock Continuous-time inverse quadratic optimal control problem.
\newblock {\em Automatica}, 117:108977, 2020.

\bibitem{li2018convex}
Yibei Li, Han Zhang, Yu~Yao, and Xiaoming Hu.
\newblock A convex optimization approach to inverse optimal control.
\newblock In {\em 2018 37th Chinese Control Conference (CCC)}, pages 257--262.
  IEEE, 2018.

\bibitem{ljung1999system}
Lennart Ljung.
\newblock {\em System Identification (2nd Ed.): Theory for the User}.
\newblock Prentice Hall PTR, USA, 1999.

\bibitem{ljung2013convexity}
Lennart Ljung and Tianshi Chen.
\newblock Convexity issues in system identification.
\newblock In {\em 2013 10th IEEE International Conference on Control and
  Automation (ICCA)}, pages 1--9, 2013.

\bibitem{lofberg2004yalmip}
Johan Lofberg.
\newblock {YALMIP}: A toolbox for modeling and optimization in {MATLAB}.
\newblock In {\em 2004 IEEE international conference on robotics and automation
  (IEEE Cat. No. 04CH37508)}, pages 284--289. IEEE, 2004.

\bibitem{menner2019constrained}
Marcel Menner, Peter Worsnop, and Melanie~N. Zeilinger.
\newblock Constrained inverse optimal control with application to a human
  manipulation task.
\newblock {\em IEEE Transactions on Control Systems Technology},
  29(2):826--834, 2021.

\bibitem{molloy2018finite}
Timothy~L Molloy, Jason~J Ford, and Tristan Perez.
\newblock Finite-horizon inverse optimal control for discrete-time nonlinear
  systems.
\newblock {\em Automatica}, 87:442--446, 2018.

\bibitem{molloy2020online}
Timothy~L Molloy, Jason~J Ford, and Tristan Perez.
\newblock Online inverse optimal control for control-constrained discrete-time
  systems on finite and infinite horizons.
\newblock {\em Automatica}, 120:109109, 2020.

\bibitem{mombaur2010human}
Katja Mombaur, Anh Truong, and Jean-Paul Laumond.
\newblock From human to humanoid locomotion -- an inverse optimal control
  approach.
\newblock {\em Autonomous robots}, 28(3):369--383, 2010.

\bibitem{pauwels2016linear}
Edouard Pauwels, Didier Henrion, and Jean-Bernard Lasserre.
\newblock Linear conic optimization for inverse optimal control.
\newblock {\em SIAM Journal on Control and Optimization}, 54(3):1798--1825,
  2016.

\bibitem{priess2014solutions}
M~Cody Priess, Richard Conway, Jongeun Choi, John~M Popovich, and Clark
  Radcliffe.
\newblock Solutions to the inverse {LQR} problem with application to biological
  systems analysis.
\newblock {\em IEEE Transactions on control systems technology},
  23(2):770--777, 2014.

\bibitem{rouot2017inverse}
J{\'e}r{\'e}my Rouot and Jean-Bernard Lasserre.
\newblock On inverse optimal control via polynomial optimization.
\newblock In {\em 2017 IEEE 56th Annual Conference on Decision and Control
  (CDC)}, pages 721--726. IEEE, 2017.

\bibitem{toumi2020tractable}
Noureddine Toumi, Roland Malham{\'e}, and Jerome Le~Ny.
\newblock A tractable mean field game model for the analysis of crowd
  evacuation dynamics.
\newblock In {\em 2020 59th IEEE Conference on Decision and Control (CDC)},
  pages 1020--1025. IEEE, 2020.

\bibitem{tsiantis2018optimality}
Nikolaos Tsiantis, Eva Balsa-Canto, and Julio~R Banga.
\newblock Optimality and identification of dynamic models in systems biology:
  an inverse optimal control framework.
\newblock {\em Bioinformatics}, 34(14):2433--2440, 2018.

\bibitem{van1998asymptotic}
Adrianus~W. van~der Vaart.
\newblock {\em Asymptotic statistics}.
\newblock Cambridge university press, Cambridge, United Kingdom, 1998.

\bibitem{westermann2020inverse}
Kevin Westermann, Jonathan Feng-Shun Lin, and Dana Kuli{\'c}.
\newblock Inverse optimal control with time-varying objectives: application to
  human jumping movement analysis.
\newblock {\em Scientific reports}, 10(1):1--15, 2020.

\bibitem{yu2021system}
Chengpu Yu, Yao Li, Hao Fang, and Jie Chen.
\newblock System identification approach for inverse optimal control of
  finite-horizon linear quadratic regulators.
\newblock {\em Automatica}, 129:109636, 2021.

\bibitem{zhang2019inverseCDC}
Han Zhang, Yibei Li, and Xiaoming Hu.
\newblock Inverse optimal control for finite-horizon discrete-time linear
  quadratic regulator under noisy output.
\newblock In {\em 2019 IEEE 58th Conference on Decision and Control (CDC)},
  pages 6663--6668. IEEE, 2019.

\bibitem{zhang2019inverse}
Han Zhang, Jack Umenberger, and Xiaoming Hu.
\newblock Inverse optimal control for discrete-time finite-horizon linear
  quadratic regulators.
\newblock {\em Automatica}, 110:108593, 2019.

\end{thebibliography}

\end{document}